\documentclass[11pt,reqno]{amsart}
\usepackage{amsmath,amssymb,mathrsfs,amsthm}
\usepackage{graphicx,cite,times,color}
\usepackage{bm}
\usepackage{subfigure}
\usepackage{framed,cleveref}
\usepackage{extarrows}
\usepackage{cases}
\newtheorem{theo}{Theorem}[section]
\newtheorem{coro}[theo]{Corollary}
\newtheorem{lemm}[theo]{Lemma}
\newtheorem{prop}[theo]{Proposition}
\newtheorem{assumption}[theo]{Assumption}
\newtheorem{rema}[theo]{Remark}

\numberwithin{equation}{section}
\crefname{prop}{Proposition}{Propositions}
\crefname{coro}{Corollary}{Corollarys}
\crefname{theo}{Theorem}{Theorems}

\newcommand{\T}{{\mathbb T}}

\newcommand{\dd}{{\rm d}}

\allowdisplaybreaks[4]
\definecolor{kkk}{RGB}{0,0,0}

\makeatletter
\def\sim@x@scale{.15}
\def\sim@y@scale{.05}
\def\sim@y@thick{.02}

\newsavebox\sim@upper
\newsavebox\sim@lower

\NewDocumentCommand{\xSim}{ O{} m }{%
	\TextOrMath{%
		\PackageError{TEST}{`\string\xSim` is valid in math mode only.}{}%
	}{
		\sbox\sim@upper{$\scriptsize #2$}
		\sbox\sim@lower{$\scriptsize #1$}
		\pgfmathparse{min(max(\wd\sim@upper/1em, \wd\sim@lower/1em, 1.0), 1.5)}
		\edef\sim@ratio{\pgfmathresult}
		\def\sim@x {\sim@x@scale * \sim@ratio}
		\def\sim@y {\sim@y@scale * \sim@ratio}
		\def\sim@@y{\sim@y@thick * \sim@ratio}
		\pgfmathparse{floor(max(\wd\sim@upper/1em, \wd\sim@lower/1em)) + 1}
		\edef\sim@wd{\pgfmathresult em}
		\mathrel{
			\begin{tikzpicture}[baseline=-.7ex]
				\filldraw[line width=.2pt] 
				(0, 0)
				.. controls +(\sim@x, \sim@y+\sim@@y) and +(-\sim@x, -\sim@y) .. 
				+(\sim@wd, 0) 
				node[midway, above] {\usebox\sim@upper} 
				node[midway, below] {\usebox\sim@lower}
				.. controls +(-\sim@x, -\sim@y-\sim@@y) and +(\sim@x, \sim@y) .. 
				(0, 0);
			\end{tikzpicture}
		}
	}
}
\makeatother

\begin{document}
	\title[Long-time weak convergence Analysis]{Long-time weak convergence Analysis of  a Semi-discrete Scheme for Stochastic Maxwell Equations}
	\maketitle
	\begin{center}
		\author{CHUCHU CHEN, JIALIN HONG, AND GE LIANG}
	\end{center}

\textbf{Abstract.} 
It is known from the monograph \cite[Chapter 5]{CHJ2023Book}  that the weak convergence analysis of numerical schemes for stochastic Maxwell equations is an unsolved problem.
This paper aims to fill the gap by establishing the long-time weak convergence analysis of the semi-implicit Euler scheme for  stochastic Maxwell equations.
Based on analyzing the regularity of transformed Kolmogorov equation associated to stochastic Maxwell equations and constructing a proper continuous adapted auxiliary process for the semi-implicit scheme, we present the long-time weak convergence analysis for this scheme and prove that the weak convergence order is one,  which is twice the strong convergence order. 
As applications of this result, we obtain the convergence order of the numerical invariant measure, the strong law of large numbers and central limit theorem related to the numerical solution, and the error estimate of the multi-level Monte Carlo estimator. As far as we know,  this is the first result on the weak convergence order for stochastic Maxwell equations.


\section{Introduction}

 Stochastic Maxwell equations model the precise microscopic origin of randomness in electromagnetic fields and have wide applications in aeronautics, electronics, biology, etc (see e.g., \cite{KS2014JMP,L2015JSP,RKT1989}).
 In the recent years, numerical approximations are proposed to discretize stochastic Maxwell equations in order to 
 study the intrinsic properties qualitatively and quantitatively.
Some papers are devoted to the constructions of structure-preserving numerical methods for stochastic Maxwell equations, for example, multi-symplectic numerical methods (cf. \cite{HJZ2014,CHZ2016}), the energy-preserving method (cf. \cite{HJZC2017}), and radial basis function collocation methods (cf. \cite{HHLS2022,H2023}).
On the convergence analysis of numerical schemes for stochastic Maxwell equations,  the existing works mainly focus on the strong convergence analysis. We refer interested readers to \cite{CHJ2019} for the semi-implicit Euler method, to \cite{CHJ2019RK} for the the Runge--Kutta methods, to \cite{CCHS2020EI} for the exponential integrators, and to \cite{Chen2021DG,SSX2022MSDG,CHJL2022ergodic,SSX2023DG} for the discontinuous Garlekin methods.
For a given numerical scheme for the stochastic Maxwell equations, both of the numerical solution and the exact solution are random variables in essence, which motivates us naturally to consider the convergence in the sense of distribution, that is, the weak convergence.  As is mentioned in \cite[Chapter 5]{CHJ2023Book}, however, there has been no work concerning the weak convergence analysis of numerical discretizations for stochastic Maxwell equations.

In this paper, we try to fill the gap and  consider the  
long-time weak convergence analysis for
the following stochastic Maxwell equations driven by multiplicative It\^o noise:
\begin{equation}\label{scv}
	\left\{
	\begin{split}
		&\dd \mathbf{E}-\nabla\times\mathbf{H}\dd t=-\sigma_0\mathbf{E}\dd t-\mathbf{J}_e(\mathbf{x},\mathbf{E},\mathbf{H})\dd t-\mathbf{J}_e^r(\mathbf{x},\mathbf{E},\mathbf{H})\dd W(t),\, &(t,\mathbf{x})\in\mathbb{R}_+\times \Theta,\\
		&\dd\mathbf{H}+\nabla\times\mathbf{E}\dd t=-\sigma_0\mathbf{H}\dd t-\mathbf{J}_m(\mathbf{x},\mathbf{E},\mathbf{H})\dd t-\mathbf{J}_m^r(\mathbf{x},\mathbf{E},\mathbf{H})\dd W(t),\, &(t,\mathbf{x})\in\mathbb{R}_+\times \Theta,\\
		&\mathbf{E}(0,\mathbf{x})=\mathbf{E}_0(\mathbf{x}),\,\,\mathbf{H}(0,\mathbf{x})=\mathbf{H}_0(\mathbf{x}),\,&\mathbf{x}\in \Theta,\\
		&\mathbf{n}\times\mathbf{E}=\mathbf{0},\,&(t,\mathbf{x})\in\mathbb{R}_+\times\partial \Theta,
	\end{split}\right.
\end{equation}
where $\Theta\subset\mathbb{R}^3$ is a bounded domain, $\mathbf{E}$ is the electric field, and $\mathbf{H}$ is the magnetic field. Here, $\{W(t)\}_{t\geq 0}$ is a $Q$-Wiener process with respect to the filtered probability space $(\Omega,\mathcal{F},\{\mathcal{F}_{t}\}_{t\geq 0},\mathbb{P})$ with $Q$ being a symmetric, positive definite operator on $L^2(\Theta)$, which characterizes the randomness from the medium. Letting $\{e_k\}_{k\in \mathbb{N}_+}$ be an orthonormal basis of the space $L^2(\Theta)$, then $W(t)$ can be represented as $W(t)=\sum_{k=1}^{\infty}Q^{\frac{1}{2}}e_k\beta_k(t),t\geq 0$, where $\{\beta_k(t)\}_{k\in\mathbb{N}_+}$ is a sequence of independent real-valued Brownian motions.
The function $\mathbf{J}:\Theta\times\mathbb{R}^3\times\mathbb{R}^3\rightarrow\mathbb{R}^3$ ($\mathbf{J}$ could be $\mathbf{J}_e,\mathbf{J}_e^r,\mathbf{J}_m$, or $\mathbf{J}_m^r$) describes a possibly nonlinear resistor, and the damping terms $\sigma_0\mathbf{E}$ and $\sigma_0\mathbf{H}$ may be induced by the conductivity of the medium or the perfectly matched layer technique. By writing \eqref{scv} into the abstract form, we obtain the uniform boundedness for the exact solution in $\mathbb{H}$ and $\mathcal{D}(M^2)$, and show the continuous dependence of the exact solution on the initial data with exponential decay rate.
As a byproduct, we prove the existence and uniqueness of the invariant measure as well as the strong law of large numbers and central limit theorem for the exact solution. 
 
 In order to inherit these properties numerically, we introduce the semi-implicit Euler scheme to discretize \eqref{scv} in the temporal direction. 
First, we obtain the uniform boundedness for the numerical solution in $\mathbb{H}$ and $\mathcal{D}(M^2)$, which combining with
 the continuous dependence of the numerical solution on the initial data,
 leads to the existence and uniqueness of invariant measure of the numerical solution.
 Further, we study the regularity of the solution $V(t,u_0)$ of  the transformed Kolmogorov equation associated to \eqref{scv},
 which shows that  the estimates of derivatives of $V(t,u_0)$ with respect to $u_0$ are of exponential growth with respect to time. 
  The key to the long-time weak convergence analysis for the semi-implicit Euler scheme lies in constructing an auxiliary adapted continuous process, which can
 balance the exponential growth of derivatives of $V$ with respect to time and avoid losing the weak convergence order in time.
 Based on the properties of the auxiliary process, we present the long-time weak convergence analysis of the semi-implicit Euler
 scheme and prove that the weak convergence order is one, which is twice the strong convergence order.
 To the best of our knowledge, this is the first result on the weak convergence analysis for numerical discretizations of stochastic Maxwell equations. 
 As applications of the weak convergence result, we obtain the convergence order of the numerical invariant measure,
 the strong law of large numbers and central limit theorem for the time-average of the numerical solution, and the error estimate of multi-level Monte Carlo estimator.
 
The rest of this paper is organized as follows. In section \ref{sdd}, we investigate the properties of the exact solution for stochastic Maxwell equations.
In section \ref{sdd2}, we introduce the semi-implicit Euler scheme and obtain the regularity, exponential stability, and the invariant measure of the numerical solution.  Section \ref{sec4} is devoted to the weak convergence analysis for the semi-implicit Euler scheme. Based on 
 the regularity of transformed Kolmogorov equation and properties of  the auxiliary process,
we show that the weak convergence order
of the scheme is one. Some applications of the weak convergence result are also given. The details of the proof to estimates of
some operators are provided in Appendix.

	\section{Properties of the solution of stochastic Maxwell equations}\label{sdd}
	Throughout this paper, we work with a real Hilbert space $\mathbb{H}:=L^2(\Theta)^3\times L^2(\Theta)^3$, endowed with the usual inner product and norm.
	Denote by $\mathbf{C}^k(\mathcal{X}_1;\mathcal{X}_2)$ (resp. $\overline{\mathbf{C}}_b^k(\mathcal{X}_1;\mathcal{X}_2)$) the space of  $k$th continuously differentiable mappings from one Hilbert space $\mathcal{X}_1$ to another Hilbert space $\mathcal{X}_2$ (resp. with bounded derivatives up to order $k$).
	 Let $\mathbf{C}_b^k(\mathcal{X}_1;\mathcal{X}_2)$
	stand for the set of bounded and $k$th continuously differentiable mappings from $\mathcal{X}_1$ to $\mathcal{X}_2$ with bounded derivatives up to order $k$. 
	When no confusion occurs, $\mathbf{C}^k(\mathcal{X}_1;\mathcal{X}_2)$ (resp. $\mathbf{C}_b^k(\mathcal{X}_1;\mathcal{X}_2)$) is simply written as $\mathbf{C}^k$ (resp. $\mathbf{C}^k_b$).
	Denote by $\mathcal{P}(\mathbb{H})$ the family of probability measures on $(\mathbb{H},\mathcal{B}(\mathbb{H}))$ and by $\nu(f):=\int_{\mathbb{H}}f\dd\nu$ the integral of the measurable functional $f$ with respect to $\nu\in\mathcal{P}(\mathbb{H})$.

Define the Maxwell operator by
\begin{equation}\label{Maxwell operator}
	M=\begin{pmatrix}
		0 & \nabla\times\\
		-\nabla\times &0
	\end{pmatrix}
\end{equation}
with the domain
\begin{equation}\label{domain}
	\begin{split}
		\mathcal{D}(M)&=\left\{\begin{pmatrix}
			\mathbf{E}\\
			\mathbf{H}
		\end{pmatrix}\in\mathbb{H}:
		M\begin{pmatrix}
			\mathbf{E}\\
			\mathbf{H}
		\end{pmatrix}=\begin{pmatrix}
			\nabla\times\mathbf{H}\\
			-\nabla\times\mathbf{E}
		\end{pmatrix}\in\mathbb{H},\,\mathbf{n}\times\mathbf{E}\Big|_{\partial \Theta}=\mathbf{0}
		\right\}\\
		&=H_{0}({\rm curl},\Theta)\times H({\rm curl},\Theta),
	\end{split}
\end{equation}
where the curl-spaces are defined by
\begin{align*}
	H({\rm curl},\Theta):=\{v\in L^2(\Theta)^3:\nabla\times v\in L^2(\Theta)^3\}
\end{align*}
and 
\begin{align*}
	H_0({\rm curl},\Theta):=\{v\in H({\rm curl},\Theta):\mathbf{n}\times v|_{\partial \Theta}=\mathbf{0}\}.
\end{align*}
The corresponding graph norm is $\|v\|_{\mathcal{D}(M)}:=\left(\|v\|_{\mathbb{H}}^2+\|Mv\|_{\mathbb{H}}^2\right)^{\frac{1}{2}}$. The Maxwell operator $M$ defined in \eqref{Maxwell operator} with domain \eqref{domain} is closed, skew-adjoint on $\mathbb{H}$, and thus generates a unitary $C_0$-group $\{e^{tM},t\in\mathbb{R}\}$ on $\mathbb{H}$ by Stone's theorem (see, for instance, \cite[Theorem B.1]{CHJ2023Book}). A frequently used property for the  Maxwell operator $M$ is $\langle Mu,u\rangle_{\mathbb{H}}=0$ for all $u\in\mathcal{D}(M)$.

Recursively, we define $\mathcal{D}(M^k)=\{u\in\mathcal{D}(M^{k-1}):M^{k-1}u\in\mathcal{D}(M)\}$ as the domain of the $k$th power of the operator $M$, $k\in\mathbb{N}$, with $\mathcal{D}(M^0)=\mathbb{H}$. The norm on $\mathcal{D}(M^k)$ is defined as
\begin{align*}
	\|v\|_{\mathcal{D}(M^k)}:=\left(\|v\|_{\mathbb{H}}^2+\|M^kv\|_{\mathbb{H}}^2\right)^{\frac{1}{2}}\qquad \forall \,v\in\mathcal{D}(M^k).
\end{align*}

Let $F:\mathbb{H}\rightarrow\mathbb{H}$ be a Nemytskij operator associated with $\mathbf{J}_e,\mathbf{J}_m$, which is defined by
\begin{align*}
	F(h)(\mathbf{x})=\begin{pmatrix}
		-\mathbf{J}_e(\mathbf{x},\mathbf{h}_1(\mathbf{x}),\mathbf{h}_2(\mathbf{x}))\\
		-\mathbf{J}_m(\mathbf{x},\mathbf{h}_1(\mathbf{x}),\mathbf{h}_2(\mathbf{x}))
	\end{pmatrix}\qquad \forall \,\mathbf{x}\in \Theta,\,h=(\mathbf{h}_1,\mathbf{h}_2)\in\mathbb{H}.
\end{align*}
For the diffusion term, we introduce the Nemytskij operator $B:\mathbb{H}\rightarrow HS(U_0;\mathbb{H})$ by
\begin{align*}
	\left(B(h)v\right)(\mathbf{x})=\begin{pmatrix}
		-\mathbf{J}_e^r(\mathbf{x},\mathbf{h}_1(\mathbf{x}),\mathbf{h}_2(\mathbf{x}))v(\mathbf{x})\\
		-\mathbf{J}_m^r(\mathbf{x},\mathbf{h}_1(\mathbf{x}),\mathbf{h}_2(\mathbf{x}))v(\mathbf{x})
	\end{pmatrix}\qquad \forall\,\mathbf{x}\in \Theta,\,h=(\mathbf{h}_1,\mathbf{h}_2)\in\mathbb{H},\,v\in U_0,
\end{align*}
where $U_0:=Q^{\frac{1}{2}}(L^2(\Theta))$. Here, $HS(\mathcal{U}_1;\mathcal{U}_2)$ denotes the separable Hilbert space of all Hilbert--Schmidt operators from one separable Hilbert space $\mathcal{U}_1$ to another separable Hilbert space $\mathcal{U}_2$, equipped with the inner product
\begin{align*}
	\langle\Gamma_1,\Gamma_2\rangle_{HS(\mathcal{U}_1;\mathcal{U}_2)}=\sum_{j=1}^{\infty}\langle\Gamma_1\eta_j,\Gamma_2\eta_j\rangle_{\mathcal{U}_2}
\end{align*}
and the corresponding norm
\begin{align*}
	\|\Gamma\|_{HS(\mathcal{U}_1;\mathcal{U}_2)}=\left(\sum_{j=1}^{\infty}\|\Gamma\eta_j\|^2_{\mathcal{U}_2}\right)^{\frac{1}{2}},
\end{align*}
where $\{\eta_j\}_{j\in\mathbb{N}_+}$ is an orthonormal basis of $\mathcal{U}_1$.

At this point, we consider the abstract form of  \eqref{scv} in the infinite-dimensional space $\mathbb{H}$:
\begin{equation}\label{www}
	\left\{\begin{split}
		&\dd u(t)=\Big(Mu(t)-\sigma_0u(t)+F(u(t))\Big)\dd t+B(u(t))\dd W(t),\\
		&u(0)=u_0.
	\end{split}\right.
\end{equation}
The mild solution of \eqref{www} is
\begin{align*}
	u(t)&=S(t)u_0+\int_{0}^{t}S(t-s)F(u(s))\dd s+\int_{0}^{t}S(t-s)B(u(s))\dd W(s),
\end{align*}
where $S(t)=e^{t(M-\sigma_0I)},t\geq 0$.

In the following, we give the assumption on the drift term $F$ and the diffusion term $B$.
\begin{assumption}\label{ass}
Let $\rho\in\mathbb{N}$. Assume that 
\begin{align}
	\langle u-v,F(u)-F(v)\rangle_{\mathcal{D}(M^\rho)}&\leq \alpha_{F,\rho}\|u-v\|^2_{\mathcal{D}(M^\rho)}\qquad\forall\, u,v\in\mathcal{D}(M^\rho),\label{ass1}\\
	\|B(u)-B(v)\|_{HS(U_0;\mathcal{D}(M^\rho))}&\leq \alpha_{B,\rho}\|u-v\|_{\mathcal{D}(M^\rho)}\qquad\forall\, u,v\in\mathcal{D}(M^\rho).\label{ass2}
\end{align}
\end{assumption}
\begin{rema}
	(i) In view of \eqref{ass1}, it is clear that 
	\begin{align}\label{cond2}
		&\langle u,F(u)\rangle_{\mathcal{D}(M^\rho)}\leq \widetilde{\alpha}_{F,\rho}(1+\|u\|^2_{\mathcal{D}(M^\rho)})\qquad\forall\, u\in\mathcal{D}(M^\rho),
	\end{align}
	where $\widetilde{\alpha}_{F,\rho}:=\frac{1}{2}\max\{1+2\alpha_{F,\rho},\|F(0)\|_{\mathcal{D}(M^\rho)}^2\}$.
	
	(ii) It follows from \eqref{ass2} that
	\begin{align}\label{cond3}
		&\|B(u)\|_{HS(U_0;\mathcal{D}(M^\rho))}\leq\widetilde{\alpha}_{B,\rho}(1+\|u\|_{\mathcal{D}(M^\rho)})\qquad\forall\, u\in\mathcal{D}(M^\rho),
	\end{align}
	where $\widetilde{\alpha}_{B,\rho}:=\max\{\alpha_{B,\rho},\|B(0)\|_{HS(U_0;\mathcal{D}(M^\rho))}\}$.
\end{rema}

	 The following proposition  gives the uniform estimate of the $r$th moment of the exact solution.
	\begin{prop}\label{smzpp}
	Let $r\geq 2$, $u_0\in L^r(\Omega;\mathbb{H})$, and Assumption 	\ref{ass} with $\rho=0$ hold. If $\sigma_0>2\widetilde{\alpha}_{F,0}+2(r-1)\widetilde{\alpha}_{B,0}^2$, then
		\begin{align}
			\sup_{t\geq 0}\mathbb{E}\big[\|u(t)\|^{r}_{\mathbb{H}}\big]\leq C\big(1+\mathbb{E}\big[\|u_0\|^{r}_{\mathbb{H}}\big]\big),
		\end{align}
	where the positive constant $C$ depends on $\sigma_0,r,\widetilde{\alpha}_{F,0}$ and $\widetilde{\alpha}_{B,0}$.
	\end{prop}
	\begin{proof}
	 Applying the  It\^o formula to $\|u(t)\|^{r}_{\mathbb{H}}$, we have
		\begin{equation}\label{smz}
		\begin{split}
			\dd\mathbb{E}\big[\|u(t)\|^{r}_{\mathbb{H}}\big]&=-r\sigma_0\mathbb{E}\big[\|u(t)\|^{r}_{\mathbb{H}}\big]\dd t+r\mathbb{E}\big[\|u(t)\|_{\mathbb{H}}^{r-2}\langle u(t),F(u(t))\rangle_{\mathbb{H}}\big]\dd t\\
			&\quad+\frac{1}{2}r(r-2)\sum_{k=1}^{\infty}\mathbb{E}\big[\|u(t)\|^{r-4}_{\mathbb{H}}|\langle u(t),B(u(t))Q^{\frac{1}{2}}e_k\rangle_{\mathbb{H}}|^2\big]\dd t\\
			&\quad+\frac{r}{2}\sum_{k=1}^{\infty}\mathbb{E}\big[\|u(t)\|^{r-2}_{\mathbb{H}}\|B(u(t))Q^{\frac{1}{2}}e_k\|^2_{\mathbb{H}}\big]\dd t,
		\end{split}
		\end{equation}
	which combining with \eqref{cond2}-\eqref{cond3} yields that
	\begin{align*}
		\dd\mathbb{E}\big[\|u(t)\|^{r}_{\mathbb{H}}\big]\leq -r\big(\sigma_0-2\widetilde{\alpha}_{F,0}-2(r-1)\widetilde{\alpha}_{B,0}^2\big)\mathbb{E}\big[\|u(t)\|^r\big]\dd t+r(\widetilde{\alpha}_{F,0}+(r-1)\widetilde{\alpha}^2_{B,0})\dd t.
	\end{align*}
	
In view of the Gronwall inequality, we finish the proof.
	\end{proof}

%\cref{smzpp}
Especially, \eqref{www} with $F(0)=0$ and $B(0)=0$ has the trivial solution, i.e., $u(t)\equiv0$ when $u(0)=0$.  Similarly to the proof of \cref{smzpp}, we can obtain the exponential stability of the moment of the exact solution for \eqref{www}.

\begin{coro}\label{stability_exact}
	Let $u_0\in  L^2(\Omega;\mathbb{H})$ and Assumption \ref{ass}  with $\rho=0$ hold. If $F(0)=0,B(0)=0,$ and $\sigma_0>\alpha_{F,0}+\frac{1}{2}\alpha^2_{B,0}$, then $\{u(t)\}_{t\geq 0}$ is exponentially stable in 
	mean square, i.e,
	\begin{align*}
		\lim_{t\rightarrow\infty}\frac{1}{t}\log\mathbb{E}\big[\|u(t)\|^2_{\mathbb{H}}\big]<-2(\sigma_0-\alpha_{F,0}-\frac{1}{2}\alpha^2_{B,0}).
	\end{align*}
\end{coro}

Next, we investigate the regularity of the exact solution in $L^2(\Omega;\mathcal{D}(M^2))$.
\begin{prop}
	Let $u_0\in L^2(\Omega;\mathcal{D}(M^2))$ and Assumption \ref{ass} with $\rho=2$ hold. If $\sigma_0>\widetilde{\alpha}_{F,2}+\widetilde{\alpha}^2_{B,2}$, then
	\begin{align*}
		\sup_{t\geq 0}\mathbb{E}\big[\|u(t)\|^2_{\mathcal{D}(M^2)}\big]\leq C(1+\mathbb{E}\big[\|u_0\|_{\mathcal{D}(M^2)}^2\big]),
	\end{align*}
where the constant $C$ depends on $\sigma_0,\widetilde{\alpha}_{F,2}$ and $\widetilde{\alpha}_{B,2}$.
\end{prop}
\begin{proof}
	We apply the It\^o formula to $\|M^2u(t)\|^2_{\mathbb{H}}$ to obtain
	\begin{align*}
		\dd\mathbb{E}\big[\|M^2u(t)\|^2_{\mathbb{H}}\big]&=-2\sigma_0\mathbb{E}\big[\|M^2u(t)\|^2_{\mathbb{H}}\big]\dd t+2\mathbb{E}\big[\langle M^2u(t),M^2F(u(t))\rangle_{\mathbb{H}}\big]\dd t\\
		&\qquad+\sum_{k=1}^{\infty}\mathbb{E}\big[\|M^2(B(u(t))Q^{\frac{1}{2}}e_k)\|^2_{\mathbb{H}}\big]\dd t.
	\end{align*}
Letting $r=2$  in \eqref{smz}, it holds that
 \begin{align*}
 	\dd\mathbb{E}\big[\|u(t)\|^2_{\mathcal{D}(M^2)}\big]&=-2\sigma_0\mathbb{E}\big[\|u(t)\|^2_{\mathcal{D}(M^2)}\big]\dd t+2\mathbb{E}\big[\langle u(t),F(u(t))\rangle_{\mathcal{D}(M^2)}\big]\dd t\\
 	&\quad+\sum_{k=1}^{\infty}\mathbb{E}\big[\|B(u(t))Q^{\frac{1}{2}}e_k\|^2_{\mathcal{D}(M^2)}\big]\dd t\\
 	&\leq -2(\sigma_0-\widetilde{\alpha}_{F,2}-\widetilde{\alpha}_{B,2}^2)\mathbb{E}\big[\|u(t)\|^2_{\mathcal{D}(M^2)}\big]\dd t+2(\widetilde{\alpha}_{F,2}+\widetilde{\alpha}_{B,2})\dd t.
 \end{align*}
Hence, we finish the proof by the Gronwall inequality.
\end{proof}

	Let 
		 $\widetilde{u}(t)$ be the solution of \eqref{www} with initial data $\widetilde{u}_0$.	Similarly to Proposition \ref{smzpp}, 
		 we obtain the continuous dependence of the solution on
		 the initial data with exponential decay rate.
	\begin{coro}\label{suu}
	Let $u_0\in L^2(\Omega;\mathbb{H})$, $\widetilde{u}_0\in L^2(\Omega;\mathbb{H})$, and Assumption 	\ref{ass}  with $\rho=0$ hold. If $\sigma_0>\alpha_{F,0}+\frac{1}{2}\alpha^2_{B,0}$, then
		\begin{align}
			\big(\mathbb{E}\big[\|u(t)-\widetilde{u}(t)\|^2_{\mathbb{H}}\big]\big)^{\frac{1}{2}}\leq e^{-(\sigma_0-\alpha_{F,0}-\frac{1}{2}\alpha^2_{B,0})t}\big(\mathbb{E}\big[\|u_0-\widetilde{u}_0\|^2_{\mathbb{H}}\big]\big)^{\frac{1}{2}}\qquad\forall\,t\in(0,\infty).
		\end{align}
	\end{coro}
	
	Let $u(t,u_0)$ be the solution of \eqref{www} at time $t$ with initial data $u_0$. Denote $P_tf(u_0):= \mathbb{E}[f(u(t,u_0))], t \geq 0$ the corresponding Markov transition semi-group. 	For $p\geq 1$ and $\gamma\in(0,1]$, 
	we define the  test functional space $\mathcal{C}_{p,\gamma}(\mathbb{H};\mathbb{R})$ the set of continuous functions on $\mathbb{H}$ endowed with the following norm
	\begin{align*}
		\|f\|_{p,\gamma}:=\sup_{h\in\mathbb{H}}\frac{|f(h)|}{1+\|h\|_{\mathbb{H}}^{p/2}}+\sup_{\substack{h_1,h_2\in \mathbb{H}\\h_1\neq h_2}}\frac{|f(h_1)-f(h_2)|}{(1\wedge\|h_1-h_2\|^\gamma_{\mathbb{H}})(1+\|h_1\|^p_{\mathbb{H}}+\|h_2\|^p_{\mathbb{H}})^{\frac{1}{2}}}.
	\end{align*}
	In view of Proposition \ref{smzpp}, Corollary \ref{suu} and \cite[Proposition 2.1]{chen2023probabilistic}, we can obtain the invariant measure, strong law of large numbers and the central limit theorem for the exact solution.
	\begin{prop}\label{drr}
		Let Assumption 	\ref{ass}  with $\rho=0$ hold.	Let  $\gamma\in(0,1]$, $r\geq 2$ and $p\geq 1$ with $2p+4\gamma\leq r$, and $u_0\in L^r(\Omega;\mathbb{H})$. If $\sigma_0>2\widetilde{\alpha}_{F,0}+2(r-1)\widetilde{\alpha}_{B,0}^2$, then $\{u(t)\}_{t\geq 0}$ admits a unique invariant measure $\pi\in\mathcal{P}(\mathbb{H})$, and fulfills the strong law of large numbers and the central limit theorem: for any $f\in \mathcal{C}_{p,\gamma}(\mathbb{H};\mathbb{R})$,
		\begin{align*}
			&\frac{1}{T}\int_{0}^{T}f(u(t))\dd t\stackrel{a.s.}{\longrightarrow}\pi(f)\qquad\text{as }T\rightarrow\infty,\\
			&\frac{1}{\sqrt{T}}\int_{0}^{T}(f(u(t))-\pi(f))\dd t\stackrel{d}{\longrightarrow}\mathcal{N}(0,v^2)\qquad\text{as }T\rightarrow\infty,
		\end{align*} 
	where $v^2:=2\pi((f-\pi(f))\int_{0}^{\infty}(P_tf-\pi(f))\dd t)$.
	\end{prop}

	\section{Temporal semi-discretization}\label{sdd2}

	In this section, we apply the semi-implicit Euler scheme to discretize \eqref{www} in the temporal direction and investigate the properties of the numerical solution. For the time interval $[0,+\infty)$, we introduce a uniform partition with step-size $\Delta t$.

Applying the semi-implicit Euler scheme to \eqref{www} in the temporal direction, we have
\begin{align}\label{semi-numerical}
	u_{k+1}-u_k=\Delta tMu_{k+1}-\sigma_0u_{k+1}\Delta t+\Delta tF(u_{k+1})+B(u_k)\Delta W_k
\end{align}
for $k\in\mathbb{N}$.
Here, the increment $\Delta W_k$ is given by
\begin{align*}
	\Delta W_k:=W(t_{k+1})-W(t_k)=\sum_{j=1}^{\infty}\left(\beta_j(t_{k+1})-\beta_j(t_k)\right)Q^{\frac{1}{2}}e_j.
\end{align*}
Let $S_{\Delta t}:=(I+\sigma_0\Delta t-M\Delta t)^{-1}$. Then \eqref{semi-numerical} is equivalent to
\begin{align}\label{ddd}
	u_{k+1}=S_{\Delta t}u_k+\Delta tS_{\Delta t}F(u_{k+1})+S_{\Delta t}B(u_k)\Delta W_k.
\end{align}

	\subsection{Properties of the numerical solution}
In this subsection, we give the uniform boundedness, the invariant measure and the exponential stability of the numerical solution $\{u_k\}_{k\geq 0}$ for  \eqref{semi-numerical}.

For $q\in\mathbb{N}_+$, denote $$\mathcal{R}_q:=\Big(\frac{q(q-1)}{2}\Big(\frac{q}{q-1}\Big)^{q-2}\Big)^{\frac{1}{2}},$$
$$\widetilde{\mathcal{R}}_q:=\big(1+4\widetilde{\alpha}_{B,0}\mathcal{R}_{2q}\big)^{q}+\sum_{i_1=0}^{q-1}\sum_{i_2=0}^{q-i_1-1}\binom{q}{i_1}\binom{q-i_1}{i_2}2^{3q-3i_1-i_2-1}\mathcal{R}_{2q}^{2q-2i_1-i_2}\widetilde{\alpha}_{B,0}^{2q-2i_1-i_2}.$$
		
				\begin{prop}\label{dfz21}
			Let $q\in\mathbb{N}_+$, $u_0\in L^{2q}(\Omega;\mathbb{H})$, Assumption \ref{ass}  with $\rho=0$ hold, and $\sigma_0>\widetilde{\alpha}_{F,0}+\frac{1}{q}\widetilde{\mathcal{R}}_q$.
			There exists a unique $\mathbb{H}$-valued $\{\mathcal{F}_{t_k}\}_{k\geq 0}$-adapted numerical solution $\{u_k\}_{k\geq 0}$ of the scheme \eqref{semi-numerical} such that for sufficiently small $\Delta t>0$,
			\begin{align}\label{dxx21}
				\sup_{ k\geq 0}\mathbb{E}\big[\|u_k\|^{2q}_{\mathbb{H}}\big]<C\big(1+\mathbb{E}\big[\|u_0\|_{\mathbb{H}}^{2q}\big]\big),
			\end{align}
		where the
		positive constant $C$ depends on $\sigma_0,q,\widetilde{\alpha}_{B,0}$ and $\widetilde{\alpha}_{F,0}$.
		\end{prop}
	\begin{proof}
		\textbf{Step 1}. First, we apply $\langle\cdot,u_{k+1}\rangle_{\mathbb{H}}$ on both sides of \eqref{semi-numerical} to obtain
		\begin{equation}\label{dxx}
			\begin{split}
			&\quad\langle u_{k+1}-u_k,u_{k+1}\rangle_{\mathbb{H}}\\
			&=-\sigma_0\Delta t\left\|u_{k+1}\right\|_{\mathbb{H}}^2+\Delta t\langle F(u_{k+1}),u_{k+1}\rangle_{\mathbb{H}}+
			\langle B(u_k)\Delta W_k,u_{k+1}-u_k\rangle_{\mathbb{H}}+
			\langle B(u_k)\Delta W_k,u_{k}\rangle_{\mathbb{H}}\\
			&\leq -\sigma_0\Delta t\left\|u_{k+1}\right\|_{\mathbb{H}}^2+\Delta t\langle F(u_{k+1}),u_{k+1}\rangle_{\mathbb{H}}\\
			&\quad+\frac{1}{2}\left\|B(u_k)\Delta W_k\right\|_{\mathbb{H}}^2+\frac{1}{2}\left\|u_{k+1}-u_k\right\|_{\mathbb{H}}^2+	\langle B(u_k)\Delta W_k,u_{k}\rangle_{\mathbb{H}},
				\end{split}
		\end{equation}
		which  combining with
		\begin{align*}
			\langle u_{k+1}-u_k,u_{k+1}\rangle_{\mathbb{H}}=\frac{1}{2}\|u_{k+1}\|_{\mathbb{H}}^2-\frac{1}{2}\|u_k\|_{\mathbb{H}}^2+\frac{1}{2}\|u_{k+1}-u_k\|^2_{\mathbb{H}}
		\end{align*}
		 leads to
		\begin{align}\label{sxx}
			\left\|u_{k+1}\right\|^2_{\mathbb{H}}&\leq \left\|u_k\right\|_{\mathbb{H}}^2-2\sigma_0\Delta t\left\|u_{k+1}\right\|_{\mathbb{H}}^2+2\Delta t\langle F(u_{k+1}),u_{k+1}\rangle_{\mathbb{H}}\notag\\
			&\quad+\left\|B(u_k)\Delta W_k\right\|_{\mathbb{H}}^2+	2\langle B(u_k)\Delta W_k,u_{k}\rangle_{\mathbb{H}}\\
			&\leq \|u_k\|^2_{\mathbb{H}}-2(\sigma_0-\widetilde{\alpha}_{F,0})\Delta t\|u_{k+1}\|_{\mathbb{H}}^2+2\widetilde{\alpha}_{F,0}\Delta t+\left\|B(u_k)\Delta W_k\right\|_{\mathbb{H}}^2+	2\langle B(u_k)\Delta W_k,u_{k}\rangle_{\mathbb{H}}.\notag
		\end{align}
	 By taking the expectation on both sides, it yields that
	\begin{align*}
		\mathbb{E}\big[\left\|u_{k+1}\right\|^2_{\mathbb{H}}\big]
		&\leq \frac{\mathbb{E}\big[\|u_k\|^2_{\mathbb{H}}\big]+2\widetilde{\alpha}_{F,0}\Delta t+\Delta t\mathbb{E}\big[\left\|B(u_k)\right\|_{HS(U_0;\mathbb{H})}^2\big]}{1+2(\sigma_0-\widetilde{\alpha}_{F,0})\Delta t}\\
		&\leq \frac{1+2\widetilde{\alpha}_{B,0}^2\Delta t}{1+2(\sigma_0-\widetilde{\alpha}_{F,0})\Delta t}\mathbb{E}\big[\|u_k\|^2_{\mathbb{H}}\big]+2(\widetilde{\alpha}_{F,0}+\widetilde{\alpha}_{B,0}^2)\Delta t.
	\end{align*}
	Thus, \eqref{dxx21} holds for $q=1$.

	\textbf{Step 2}. By induction, we assume that
	\begin{align}\label{szu}
			\sup_{ k\geq 0}\mathbb{E}\left[\|u_k\|^{2(q-1)}_{\mathbb{H}}\right]<C(1+\mathbb{E}\big[\|u_0\|_{\mathbb{H}}^{2(q-1)}\big]).
	\end{align}
In view of \eqref{sxx} and $(1+2(\sigma_0-\widetilde{\alpha}_{F,0})\Delta t)^q\geq 1+2q(\sigma_0-\widetilde{\alpha}_{F,0})\Delta t$, we have
	\begin{equation}\label{khd}
		\begin{split}
		\mathbb{E}\big[\|u_{k+1}\|^{2q}_{\mathbb{H}}\big]&\leq \frac{\mathbb{E}\Big[\Big( \|u_k\|^2_{\mathbb{H}}+\big(2\widetilde{\alpha}_{F,0}\Delta t+\|B(u_k)\Delta W_k\|^2_{\mathbb{H}}\big)+2\langle B(u_k)\Delta W_k,u_{k}\rangle_{\mathbb{H}}\Big)^q\Big]}{1+2q(\sigma_0-\widetilde{\alpha}_{F,0})\Delta t}\\
		&=\frac{\mathbb{E}\big[\|u_k\|_{\mathbb{H}}^{2q}\big]+\sum_{i_1=0}^{q-1}\sum_{i_2=0}^{q-i_1-1}\binom{q}{i_1}\binom{q-i_1}{i_2}2^{i_2}\mathbb{E}\big[\mathcal{S}_{k,i_1,i_2}\big]+\sum_{i=0}^{q-1}\binom{q}{i}2^{q-i}\mathbb{E}\big[\widetilde{\mathcal{S}}_{k,i}\big]}{1+2q(\sigma_0-\widetilde{\alpha}_{F,0})\Delta t},
			\end{split}
	\end{equation}
%	Notice that
%	\begin{align*}
%		&\quad\mathbb{E}\Big[\Big( \|u_k\|^2_{\mathbb{H}}+\big(2\widetilde{\alpha}_{F,0}\Delta t+\|B(u_k)\Delta W_k\|^2_{\mathbb{H}}\big)+2\langle B(u_k)\Delta W_k,u_{k}\rangle_{\mathbb{H}}\Big)^q\Big]\\
%		&=\mathbb{E}\big[\|u_k\|_{\mathbb{H}}^{2q}\big]+\sum_{i_1=0}^{q-1}\sum_{i_2=0}^{q-i_1-1}\binom{q}{i_1}\binom{q-i_1}{i_2}2^{i_2}\mathbb{E}\big[\mathcal{S}_{k,i_1,i_2}\big]+\sum_{i=0}^{q-1}\binom{q}{i}2^{q-i}\mathbb{E}\big[\widetilde{\mathcal{S}}_{k,i}\big],
%	\end{align*}
where
\begin{align*}
	\mathcal{S}_{k,i_1,i_2}&:=\|u_k\|^{2i_1}_{\mathbb{H}}\langle B(u_k)\Delta W_k,u_k\rangle_{\mathbb{H}}^{i_2}\big(2\widetilde{\alpha}_{F,0}\Delta t+\|B(u_k)\Delta W_k\|^2_{\mathbb{H}}\big)^{q-i_1-i_2},\\
	\widetilde{\mathcal{S}}_{k,i}&:=\|u_k\|^{2i}_{\mathbb{H}}\langle B(u_k)\Delta W_k,u_k\rangle_{\mathbb{H}}^{q-i}.
\end{align*}

	Notice that for $i=0,\cdots,q-2$, it follows from the Burkholder--Davis--Gundy-type inequality  that
	\begin{align*}
		\big|\mathbb{E}\big[\widetilde{\mathcal{S}}_{k,i}\big]\big|&\leq \big(\mathbb{E}\big[\|u_k\|_{\mathbb{H}}^{2q}\big]\big)^{\frac{q+i}{2q}}\big(\mathbb{E}\big[\|B(u_k)\Delta W_k\|_{\mathbb{H}}^{2q}\big]\big)^{\frac{q-i}{2q}}
		\leq\big(2\widetilde{\alpha}_{B,0}\mathcal{R}_{2q}\big)^{q-i} \big(1+\mathbb{E}\big[\|u_k\|_{\mathbb{H}}^{2q}\big]\big)\Delta t,
	\end{align*}
	which combining with $\mathbb{E}\big[\widetilde{\mathcal{S}}_{k,q-1}\big]=0$ leads to
\begin{align}\label{swq}
	\Big|\mathbb{E}\Big[\sum_{i=0}^{q-1}\binom{q}{i}2^{q-i}\widetilde{\mathcal{S}}_{k,i}\Big]\Big|\leq \big(1+4\widetilde{\alpha}_{B,0}\mathcal{R}_{2q}\big)^{q}\Delta t+\big(1+4\widetilde{\alpha}_{B,0}\mathcal{R}_{2q}\big)^{q}\Delta t\mathbb{E}\big[\|u_k\|_{\mathbb{H}}^{2q}\big].
\end{align}
	
Besides, one obtains that  for $i_1=0,\cdots,q-1,i_2=0,\cdots,q-i_1-1$,
	\begin{align*}
	\big|\mathbb{E}\big[S_{k,i_1,i_2}\big]\big|&\leq 2^{2q-2i_1-2i_2-1}\widetilde{\alpha}_{F,0}^{q-i_1-i_2}\Delta t\mathbb{E}\big[\|u_k\|_{\mathbb{H}}^{2i_1+i_2}\|B(u_k)\Delta W_k\|_{\mathbb{H}}^{i_2}\big]\\
	&\quad+2^{q-i_1-i_2-1}\mathbb{E}\big[\|u_k\|^{2i_1+i_2}_{\mathbb{H}}\|B(u_k)\Delta W_k\|_{\mathbb{H}}^{2q-2i_1-i_2}\big]\\
	&\leq2^{2q-2i_1-i_2-1} \mathcal{R}_{2q-2}^{i_2}\widetilde{\alpha}_{B,0}^{i_2}\widetilde{\alpha}_{F,0}^{q-i_1-i_2}\Delta t\big(1+\mathbb{E}\big[\|u_k\|^{2q-2}_{\mathbb{H}}\big]\big)\\
	&\quad+2^{3q-3i_1-2i_2-1}\mathcal{R}_{2q}^{2q-2i_1-i_2}\widetilde{\alpha}_{B,0}^{2q-2i_1-i_2}\Delta t\Big(1+\mathbb{E}\big[\|u_k\|^{2q}_{\mathbb{H}}\big]\Big)\\
	&\leq 2^{3q-3i_1-2i_2-1}\mathcal{R}_{2q}^{2q-2i_1-i_2}\widetilde{\alpha}_{B,0}^{2q-2i_1-i_2}\Delta t\mathbb{E}\big[\|u_k\|^{2q}_{\mathbb{H}}\big]+C\Delta t,
	\end{align*}
which implies that
\begin{equation}\label{swq2}
	\begin{split}
	&\quad\Big|\mathbb{E}\Big[\sum_{i_1=0}^{q-1}\sum_{i_2=0}^{q-i_1-1}\binom{q}{i_1}\binom{q-i_1}{i_2}2^{i_2}\mathcal{S}_{k,i_1,i_2}\Big]\Big|\\
	&\leq \sum_{i_1=0}^{q-1}\sum_{i_2=0}^{q-i_1-1}\binom{q}{i_1}\binom{q-i_1}{i_2}2^{3q-3i_1-i_2-1}\mathcal{R}_{2q}^{2q-2i_1-i_2}\widetilde{\alpha}_{B,0}^{2q-2i_1-i_2}\Delta t\mathbb{E}\big[\|u_k\|^{2q}_{\mathbb{H}}\big]+C\Delta t.
		\end{split}
\end{equation}

By \eqref{swq} and \eqref{swq2}, we arrive at
\begin{align*}
	\mathbb{E}\big[\|u_{k+1}\|^{2q}_{\mathbb{H}}\big]&\leq \frac{1+\widetilde{\mathcal{R}}_q\Delta t}{1+2q(\sigma_0-\widetilde{\alpha}_{F,0})\Delta t}\mathbb{E}\big[\|u_k\|_{\mathbb{H}}^{2q}\big]+C\Delta t\\
	&\leq e^{-q(\sigma_0-\widetilde{\alpha}_{F,0}-\frac{1}{q}\widetilde{\mathcal{R}}_q)\Delta t}\mathbb{E}\big[\|u_k\|_{\mathbb{H}}^{2q}\big]+C\Delta t,
\end{align*}
	where $\Delta t\in(0, \frac{1}{q(\sigma_0-\widetilde{\alpha}_{F,0})})$.

Thus, \eqref{dxx21} is proved by the discrete Gronwall inequality and the  induction argument.
	\end{proof}

		Similarly to Corollary \ref{stability_exact},  the  solution of the semi-implicit Euler scheme \eqref{semi-numerical} is also exponentially stable in mean square.
		
		\begin{coro}\label{stability_numerical}
			Let $u_0\in L^2(\Omega;\mathbb{H})$ and Assumption \ref{ass} with $\rho=0$ hold. If $F(0)=0$, $B(0)=0$ and $\sigma_0>\alpha_{F,0}+\alpha^2_{B,0}$, then for sufficiently small $\Delta t>0$, $\{u_k\}_{k\geq 0}$ is exponentially stable in 
			mean square, i.e,
			\begin{align*}
				\lim_{k\rightarrow\infty}\frac{1}{k\Delta t}\log\mathbb{E}\big[\|u_k\|^2_{\mathbb{H}}\big]<-(\sigma_0-\alpha_{F,0}-\alpha^2_{B,0}).
			\end{align*}
		\end{coro}

Let $\{\widetilde{u}_k\}_{k\geq 0}$ be the solution of \eqref{semi-numerical} with initial data $\widetilde{u}_0$. Similarly to \cref{suu},  we can show the continuous dependence of the numerical solution on
the initial data with exponential decay rate.
	\begin{coro}\label{hf}
				Let $u_0\in L^2(\Omega;\mathbb{H})$,$\widetilde{u}_0\in L^2(\Omega;\mathbb{H})$, and Assumption \ref{ass}  with $\rho=0$ hold.  If $\sigma_0>\alpha_{F,0}+\alpha_{B,0}^2$, then
		\begin{align}
			\mathbb{E}\big[\|u_k-\widetilde{u}_k\|^2_{\mathbb{H}}\big]\leq e^{-(\sigma_0-\alpha_{F,0}-\alpha_{B,0}^2)k\Delta t}\mathbb{E}\big[\|u_0-\widetilde{u}_0\|^2_{\mathbb{H}}\big]\qquad\forall\,k\in\mathbb{N}_+.
		\end{align}
	\end{coro}
	In view of  \cref{dfz21},  \cref{hf} and \cite[Proposition 4.1]{chen2023probabilistic}, we can obtain the invariant measure for the numerical solution.
\begin{prop}
	Let  $u_0\in L^2(\Omega;\mathbb{H})$, Assumption \ref{ass} with $\rho=0$ hold, and $\sigma_0>\widetilde{\alpha}_{F,0}+2\widetilde{\alpha}_{B,0}^2$. For sufficiently small $\Delta t>0$, $\{u_k\}_{k\in\mathbb{N}}$ admits a unique invariant measure $\pi^{\Delta t}\in\mathcal{P}(\mathbb{H})$.
\end{prop}

%	Denote $\widetilde{C}_{2+\gamma}$ the positive constants satisfying
%\begin{align*}
%	&\|fg\|_{\mathcal{D}(M^2)}\leq \widetilde{C}_{2+\gamma}\|f\|_{H^{2+\gamma}(D)}\|g\|_{\mathcal{D}(M^2)}\qquad\forall\, f\in H^\gamma(D),g\in\mathcal{D}(M^2).
%\end{align*}
%	Let $\widehat{K}^{(2)}:=C_F^2+2\widetilde{C}^2_{2+\gamma}C_B^2\|Q^{\frac{1}{2}}\|^2_{HS(L^2(D);H^{2+\gamma}(D))}+\frac{1}{2}.$
		Similarly to \cref{dfz21}, we can prove uniform boundedness of numerical solution in $L^2(\Omega;\mathcal{D}(M^2))$.
		\begin{prop}\label{dfz3}
Let $u_0\in L^2(\Omega;\mathcal{D}(M^2))$, Assumption \ref{ass}  with $\rho=2$ hold, and
 $\sigma_0>\widetilde{\alpha}_{F,2}+2\widetilde{\alpha}_{B,2}^2$.
 For sufficiently small $\Delta t>0$,
	\begin{align}\label{dxx2}
		\sup_{ k\geq 0}\mathbb{E}\big[\|u_k\|^2_{\mathcal{D}(M^2)}\big]\leq C\big(1+\mathbb{E}\big[\|u_0\|^2_{\mathcal{D}(M^2)}\big]\big),
	\end{align}
where the positive constant  $C$ depends on $\sigma_0,\widetilde{\alpha}_{F,2}$, and $\widetilde{\alpha}_{B,2}$.
\end{prop}
\begin{proof}
	We apply $M^2$ and $\langle\cdot,M^2u_{k+1}\rangle_{\mathbb{H}}$ on both sides of \eqref{semi-numerical} to obtain that
	\begin{align*}
	&\quad\langle M^2u_{k+1}-M^2u_k,M^2u_{k+1}\rangle_{\mathbb{H}}\\
		&=-\sigma_0\Delta t\|M^2u_{k+1}\|^2_{\mathbb{H}}+\Delta t\langle M^2F(u_{k+1}),M^2u_{k+1}\rangle_{\mathbb{H}}\\
		&\quad+\langle M^2(B(u_k)\Delta W_k),M^2u_{k+1}-M^2u_k\rangle_{\mathbb{H}}+\langle M^2(B(u_k)\Delta W_k),M^2u_k\rangle_{\mathbb{H}}\\
		&\leq -\sigma_0\Delta t\|M^2u_{k+1}\|_{\mathbb{H}}^2+\Delta t\langle M^2F(u_{k+1}),M^2u_{k+1}\rangle_{\mathbb{H}}\\
		&\quad+\frac{1}{2}\|M^2(B(u_k)\Delta W_k)\|^2_{\mathbb{H}}+\frac{1}{2}\|M^2u_{k+1}-M^2u_k\|_{\mathbb{H}}^2+\langle M^2(B(u_k)\Delta W_k),M^2u_k\rangle_{\mathbb{H}},
	\end{align*}
which combining with \eqref{dxx} leads to
\begin{align*}
\mathbb{E}\big[\|u_{k+1}\|^2_{\mathcal{D}(M^2)}\big]&\leq \mathbb{E}\big[\|u_k\|^2_{\mathcal{D}(M^2)}\big]-2\sigma_0\Delta t\mathbb{E}\big[\|u_{k+1}\|^2_{\mathcal{D}(M^2)}\big]\\
&\quad+2\Delta t\mathbb{E}\big[\langle F(u_{k+1}),u_{k+1}\rangle_{\mathcal{D}(M^2)}\big]+\mathbb{E}\big[\|B(u_k)\|^2_{HS(U_0;\mathcal{D}(M^2))}\big]\Delta t\\
&\leq (1+2\widetilde{\alpha}_{B,2}^2\Delta t)\mathbb{E}\big[\|u_k\|^2_{\mathcal{D}(M^2)}\big]-2(\sigma_0-\widetilde{\alpha}_{F,2})\Delta t\mathbb{E}\big[\|u_{k+1}\|^2_{\mathcal{D}(M^2)}\big]+2(\widetilde{\alpha}_{F,2}+\widetilde{\alpha}^2_{B,2})\Delta t.
\end{align*}

	Hence, it holds that
	\begin{align*}
		\mathbb{E}\big[\|u_{k+1}\|^2_{\mathcal{D}(M^2)}\big]&\leq \frac{1+2\widetilde{\alpha}^2_{B,2}\Delta t}{1+2(\sigma_0-\widetilde{\alpha}_{F,2})\Delta t}\mathbb{E}\big[\|u_k\|^2_{\mathcal{D}(M^2)}\big]+C\Delta t\\
		&\leq e^{-(\sigma_0-\widetilde{\alpha}_{F,2}-2\widetilde{\alpha}^2_{B,2})\Delta t}\mathbb{E}\big[\|u_k\|^2_{\mathcal{D}(M^2)}\big]+C\Delta t,
	\end{align*}
where $\Delta t\in(0, \frac{1}{\sigma_0-\widetilde{\alpha}_{F,2}})$. Thus we finish the proof by the discrete Gronwall inequality.
\end{proof}

	\section{Long-time weak convergence analysis}\label{sec4}
	In this section, we first study the regularity of the Kolmogorov equation of \eqref{www}, which combining with some estimates for operators establishes the long-time weak convergence analysis for the semi-implicit Euler scheme \eqref{semi-numerical}. It is shown that the weak convergence order is twice the strong convergence order. This is the first result about the weak convergence order of the numerical scheme for stochastic Maxwell equations.
	As the applications of the weak convergence result, we obtain the error estimate between the invariant measures $\pi$ and $\pi^{\Delta t}$, the strong law of large numbers and central limit theorem for the time average of the numerical solution, and the convergence order of multi-level Monte Carlo estimator.
	
		\subsection{Kolmogorov equation}
	In this part, we focus on  the regularity of the Kolmogorov equation of \eqref{www}. 
	Let $u(t,u_0)$ be the exact solution of \eqref{www} at time $t$ with initial data $u_0$. Denote  $U(t,u_0):=\mathbb{E}\big[\varphi\big(u(t,u_0)\big)\big]$ with $\varphi\in \mathbf{C}_b^3(\mathbb{H};\mathbb{R})$. It is well known that $U(t,u_0)$ satisfies the following Kolmogorov equation (see, e.g., \cite[Chapter 9]{PZ2014}):
	\begin{align*}
		\frac{\dd U}{\dd t}(t,u)&=\langle DU(t,u),Mu-\sigma_0u+F(u)\rangle
		+\frac{1}{2}{\rm Tr}[D^2U(t,u)(B(u)Q^{\frac{1}{2}})(B(u)Q^{\frac{1}{2}})^*],
	\end{align*}
	where $U(0,u)=\varphi(u)$, $DU(t,u)$  and $D^2U(t,u)$ denote the first and second order Fr\'echet derivatives of $U(t,u)$ with respect to $u$ respectively. To eliminate the operator $M-\sigma_0I$, we consider $V(t,v):=U(t,S(-t)v)$, which satisfies
	\begin{equation}\label{spp1}
		\begin{split}
			\frac{\dd V}{\dd t}(t,v)
			&=\langle DV(t,v),\Upsilon(t,v)\rangle+\frac{1}{2}{\rm Tr}[D^2V(t,v)\Lambda(t,v)\Lambda(t,v)^*],
		\end{split}
	\end{equation}
where $\Upsilon(t,v):=S(t)F(S(-t)v)$ and $\Lambda(t,v):=S(t)B(S(-t)v)Q^{\frac{1}{2}}$.

	\begin{assumption}\label{asss}
		Assume that $F\in\overline{\mathbf{C}}_b^3(\mathbb{H};\mathbb{H})$ and $B\in\overline{\mathbf{C}}_b^3(\mathbb{H};HS(U_0;\mathbb{H}))$.
	\end{assumption}
	
	Let
	\begin{align*}
		K^{(1)}_p:=\|DF\|_{\mathbf{C}_b}+\frac{2p-1}{2}\|DB\|^2_{\mathbf{C}_b}\qquad\forall \,p\in\{1,2,3,4\},
	\end{align*}
	\begin{align*}
		K_p^{(2)}:=\|DF\|_{\mathbf{C}_b}+\frac{p+1}{4}\|D^2F\|_{\mathbf{C}_b}+
		(2p-1)\|DB\|^2_{\mathbf{C}_b}+\frac{3(p-1)}{2}\|D^2B\|^2_{\mathbf{C}_b}\qquad\forall\,p\in\{1,2\},
	\end{align*}
	\begin{align*}
		K^{(3)}:=\frac{1}{2}\|D^3F\|_{\mathbf{C}_b}+\|DF\|_{\mathbf{C}_b}+\frac{3}{2}\|D^2F\|_{\mathbf{C}_b}+\frac{3}{2}\|DB\|^2_{\mathbf{C}_b},\qquad K^{(4)}:=\max\{K_3^{(1)},K_1^{(2)}, K^{(3)}\}.
	\end{align*}

	The following lemma gives the regularity of functions $U$ and $V$.

	\begin{lemm}
		Let Assumption \ref{asss} hold and $\sigma_0>\widehat{K}:=\max\big\{K_4^{(1)},K_2^{(2)},K^{(3)},2K_4^{(1)}-K_2^{(2)}\big\}$. The functions $U$ and $V$ are continuous in time. For any $t\in (0,\infty)$, $U(t,\cdot)$ and $V(t,\cdot)$ take values in $\mathbf{C}^3_b(\mathbb{H};\mathbb{R})$.
Moreover, we have the following estimates:\\
$(i)$	$
	\|DU(t,u_0)\|_{\mathcal{L}(\mathbb{H};\mathbb{R})}\leq \|D\varphi\|_{\mathbf{C}_b}e^{-(\sigma_0-K_1^{(1)})t}$ and $	\|DV(t,v)\|_{\mathcal{L}(\mathbb{H};\mathbb{R})}\leq \|D\varphi\|_{\mathbf{C}_b}e^{K^{(1)}_1t}.$\\
	$(ii)$ $\|D^2U(t,u_0)\|_{\mathcal{L}(\mathbb{H}\times\mathbb{H};\mathbb{R})}\leq Ce^{-(\sigma_0-K^{(2)}_1)t}$ and $	\|D^2V(t,u_0)\|_{\mathcal{L}(\mathbb{H}\times\mathbb{H};\mathbb{R})}\leq Ce^{(\sigma_0+K^{(2)}_1)t}.$\\
	$(iii)$ $	\|D^3U(t,u_0)\|_{\mathcal{L}(\mathbb{H}\times\mathbb{H}\times\mathbb{H};\mathbb{R})}\leq Ce^{-(\sigma_0-K^{(4)})t}$ and $
	\|D^3V(t,u_0)\|_{\mathcal{L}(\mathbb{H}\times\mathbb{H}\times\mathbb{H};\mathbb{R})}\leq Ce^{(2\sigma_0+K^{(4)})t}$.
	\end{lemm}
	\begin{proof}
		
	(i) For every $h\in\mathbb{H}$, we have
		\begin{align*}
			DU(t,u_0)h=\mathbb{E}\big[D\varphi\big(u(t,u_0)\big)\eta^h(t)\big]\qquad \forall~t\in(0,\infty),
		\end{align*}
		where
		\begin{equation}
			\left\{\begin{split}
				\dd \eta^h(t)&=\left(M\eta^h(t)-\sigma_0\eta^h(t)+DF(u(t,u_0))\eta^h(t)\right)\dd t+DB(u(t,u_0))\eta^h(t)\dd W(t),\\
				\eta^h(0)&=h.
			\end{split}\right.
		\end{equation}
		
		By applying the It\^o formula to $\|\eta^h\|_{\mathbb{H}}^{2p}$ with $1\leq p\leq 4$ and taking the expectation, it holds that
		\begin{align*}
			\dd \mathbb{E}\big[\|\eta^h(t)\|_{\mathbb{H}}^{2p}\big]
			&=-2p\sigma_0\mathbb{E}\big[\|\eta^h(t)\|_{\mathbb{H}}^{2p}\big]\dd t
			+2p\mathbb{E}\big[\|\eta^h(t)\|_{\mathbb{H}}^{2p-2}\langle\eta^h(t),DF(u(t,u_0))\eta^h(t)\rangle_{\mathbb{H}}\big]\dd t
			\\
			&\quad+p\sum_{k=1}^{\infty}\mathbb{E}\big[\|\eta^h(t)\|^{2p-2}_{\mathbb{H}}\|DB(u(t,u_0))\eta^h(t)Q^{\frac{1}{2}}e_k\|_{\mathbb{H}}^2\big]\dd t\\
			&\quad+2p(p-1)\sum_{k=1}^{\infty}\mathbb{E}\big[\|\eta^h(t)\|^{2p-4}_{\mathbb{H}}\big|\langle\eta^h(t),DB(u(t,u_0))\eta^h(t)Q^{\frac{1}{2}}e_k\rangle_{\mathbb{H}}\big|^2\big]\dd t\\
			&\leq-2p(\sigma_0-K_p^{(1)})\mathbb{E}\big[\|\eta^h(t)\|^{2p}_{\mathbb{H}}\big]\dd t.
		\end{align*}
		In view of the Gronwall inequality, one has
		\begin{align}\label{dmm1}
			\mathbb{E}\big[\|\eta^h(t)\|_{\mathbb{H}}^{2p}\big]\leq \|h\|^{2p}_{\mathbb{H}}e^{-2p(\sigma_0-K_{p}^{(1)})t}\qquad\forall~t\in(0,\infty).
		\end{align}

		Hence, 
		\begin{align*}
			|DU(t,u_0)(h)|
			&\leq\|D\varphi\|_{\mathbf{C}_b}\big(\mathbb{E}\big[\|\eta^h(t)\|^2_{\mathbb{H}}\big]\big)^{\frac{1}{2}}
			\leq \|D\varphi\|_{\mathbf{C}_b}\|h\|_{\mathbb{H}}e^{-(\sigma_0-K_1^{(1)})t}.
		\end{align*}
		
		Since $DV(t,v)h=DU(t,S(-t)v)S(-t)h$, 
		we have
		\begin{align*}
			|DV(t,v)(h)|
			&\leq\|DU(t,S(-t)v)\|_{\mathcal{L}(\mathbb{H};\mathbb{R})}\|S(-t)\|_{\mathcal{L}(\mathbb{H};\mathbb{H})}\|h\|_{\mathbb{H}}
			\leq e^{K_1^{(1)}t} \|D\varphi\|_{\mathbf{C}_b}\|h\|_{\mathbb{H}}.
		\end{align*}

		(ii)  For every $h\in\mathbb{H}$, notice that
		\begin{align*}
			D^2U(t,u_0)(h,h)=\mathbb{E}\big[D^2\varphi( u(t,u_0))\big(\eta^h(t),\eta^h(t)\big)+D\varphi(u(t,u_0))\xi^h(t)\big]\qquad\forall~t\in(0,\infty),
		\end{align*}
		where
		\begin{equation}
			\left\{\begin{split}
				\dd \xi^h(t)&=\Big(M\xi^h(t)-\sigma_0\xi^h(t)+D^2F(u(t,u_0))(\eta^h(t),\eta^h(t))+DF(u(t,u_0))\xi^h(t)\Big)\dd t\\
				&\qquad\quad\,+\Big(D^2B(u(t,u_0))(\eta^h(t),\eta^h(t))+DB(u(t,u_0))\xi^h(t)\Big)\dd W(t),\\
				\xi^h(0)&=0.
			\end{split}\right.
		\end{equation}
		
		We apply the It\^o formula to $\|\xi^h\|_{\mathbb{H}}^2$ and take the expectation to obtain that
		\begin{align*}
			\dd \mathbb{E}\left[\|\xi^h(t)\|_{\mathbb{H}}^2\right]&=2\mathbb{E}\left[\left\langle\xi^h(t),-\sigma_0\xi^h(t)+D^2F(u(t,u_0))\big(\eta^h(t),\eta^h(t)\big)+DF(u(t,u_0))\xi^h(t)\right\rangle_{\mathbb{H}}\right]\dd t\\
			&\quad+\sum_{k=1}^{\infty}\mathbb{E}\left[\left\|\left(D^2B(u(t,u_0))\big(\eta^h(t),\eta^h(t)\big)+DB(u(t,u_0))\xi^h(t)\right)Q^{\frac{1}{2}}e_k\right\|^2_{\mathbb{H}}\right]\dd t\\
			&\leq -2(\sigma_0-K_1^{(2)})\mathbb{E}\left[\|\xi^h(t)\|^2_{\mathbb{H}}\right]\dd t+C\mathbb{E}\big[\|\eta^h(t)\|^4_{\mathbb{H}}\big]\dd t\\
			&\leq -2(\sigma_0-K_1^{(2)})\mathbb{E}\left[\|\xi^h(t)\|^2_{\mathbb{H}}\right]\dd t+C\|h\|^4_{\mathbb{H}}e^{-4(\sigma_0-K_2^{(1)})t}\dd t.
		\end{align*}
		By the Gronwall inequality, we have
		\begin{align*}
			\mathbb{E}\left[\|\xi^h(t)\|_{\mathbb{H}}^2\right]&\leq C\|h\|^4_{\mathbb{H}}\int_{0}^{t}e^{-2(\sigma_0-K_1^{(2)})(t-s)}e^{-4(\sigma_0-K_2^{(1)})s}\dd s
			\leq C\|h\|^4_{\mathbb{H}}e^{-2(\sigma_0-K_1^{(2)})t},
		\end{align*}
		 which implies that
		\begin{align*}
			|D^2U(t,u_0)(h,h)|&\leq \|D^2\varphi\|_{\mathbf{C}_b}\mathbb{E}\big[\|\eta^h(t)\|^2_{\mathbb{H}}\big]+\|D\varphi\|_{\mathbf{C}_b}\mathbb{E}\big[\|\xi^h(t)\|_{\mathbb{H}}\big]
			\leq C\|h\|^2_{\mathbb{H}}e^{-(\sigma_0-K^{(2)}_1)t}.
		\end{align*}
		
		Notice that
		\begin{align*}
			D^2V(t,v)(h,h)=D^2U(t,S(-t)v)(S(-t)h,S(-t)h),
		\end{align*}
		which implies that
		\begin{align*}
			|D^2V(t,v)(h,h)|
			&\leq \|D^2U(t,S(-t)u_0)\|_{\mathcal{L}(\mathbb{H}\times\mathbb{H};\mathbb{R})}\|S(-t)\|^2_{\mathcal{L}(\mathbb{H};\mathbb{H})}\|h\|_{\mathbb{H}}^2\leq C\|h\|_{\mathbb{H}}^2e^{(\sigma_0+K^{(2)}_1)t}.
		\end{align*}
	
	Since the proof of (iii) is similar to (i) and (ii), we omit it here. Thus we finish the proof.
	\end{proof}

	\subsection{Weak convergence order}
	
In order to derive the error estimate of \eqref{semi-numerical}, we need to introduce the following lemma, whose proof is given in Appendix A.
\begin{lemm}\label{sfu4}
	Let $g\in(0,2)$ and $\widetilde{g}\in(0,1)$.\\
	$(1)$  There exists a constant $\Delta t_g>0$ such that for every $\Delta t\in(0,\Delta t_g)$,
	  $$\|S^k_{\Delta t}\|_{\mathcal{L}(\mathcal{D}(M^i);\mathcal{D}(M^i))}\leq e^{-\frac{g}{2}k\sigma_0\Delta t}\qquad\forall\,i\in\{0,1,2\}.$$
	$(2)$ There exists  a constant $C>0$ such that  $$\|S(t)-I\|_{\mathcal{L}(\mathcal{D}(M);\mathbb{H})}\leq Ct\qquad\forall\,t\in(0,\infty).$$
	$(3)$ There exists  a constant $\Delta t_{\widetilde{g}}>0$ such that for every $\Delta t\in(0,\Delta t_{\widetilde{g}})$, 
	\begin{align}\label{skk4}
		\|S_{\Delta t}^{k}-S(t)\|_{\mathcal{L}(\mathcal{D}(M^2);\mathbb{H})}\leq C\Delta te^{-\frac{\widetilde{g}}{2}k\sigma_0\Delta t}\qquad\forall\, t\in[t_{(k-1)\vee0},t_{k+1}].
	\end{align}
\end{lemm}

\begin{assumption}\label{assss}
	Assume that for every $h\in\mathcal{D}(M^2)$,
	\begin{align*}
		\|F(h)\|_{\mathcal{D}(M^2)}&\leq L_F(1+\|h\|_{\mathcal{D}(M^2)}).
		%			\|B(u)\|_{\mathcal{D}(M^2)}&\leq L_B(1+\|u\|_{\mathcal{D}(M^2)}).
	\end{align*}
\end{assumption}
\begin{theo}\label{zx}
	Let Assumption \ref{ass} with $\rho\in\{0,2\}$, Assumption \ref{asss} and Assumption \ref{assss} hold.
 Let $u_0\in L^2(\Omega;\mathcal{D}(M^2))\cap L^4(\Omega;\mathbb{H})$, $g\in(1,2)$, $\widetilde{g}\in(0,1)$ and $\sigma_0>0$ satisfying  
 $$\widetilde{g}+2g-4>0\quad\text{and}\quad \sigma_0>\max\Big\{\frac{2K^{(4)}}{3g-4},\frac{2K^{(2)}_1}{\widetilde{g}+2g-4},\widehat{K},\widetilde{\alpha}_{F,2}+2\widetilde{\alpha}_{B,2}^2,2\widetilde{\alpha}_{F,0}+\frac{1}{2}\widetilde{\mathcal{R}}_2\Big\}.$$
 For any $\varphi\in \mathbf{C}_b^3(\mathbb{H};\mathbb{R})$, there exists a positive constant $C$ such that for sufficiently small $\Delta t>0$,
	\begin{align*}
		\sup_{N\geq 0}\big|\mathbb{E}\big[\varphi(u_N)\big]-\mathbb{E}\big[\varphi(u(t_N))\big]\big|\leq C\Delta t.
	\end{align*}
\end{theo}
\begin{proof}
	\textbf{Step 1}. 
 We define the continuous process $\{\widetilde{v}^{[k]}(t)\}_{t\in[t_k,t_{k+1}]}$ by
	\begin{equation}\label{dddwww4}
		\left\{\begin{split}
			\dd \widetilde{v}^{[k]}(t)&=S_{\Delta t}^{N-k}F(S_{\Delta t}u_k)\dd t+J_k\dd W(t),\quad t\in[t_k,t_{k+1}],\\
			\widetilde{v}^{[k]}(t_k)&=v_k,
		\end{split}\right.
	\end{equation}
where $J_k:=S_{\Delta t}^{N-k-1}(I+\Delta tS_{\Delta t}DF(S_{\Delta t}u_k))S_{\Delta t}B(u_k)$ and $	v_k:=S_{\Delta t}^{N-k}u_k$.  

\textbf{(1a). Estimate of }$	\mathbb{E}\big[\|v_{k+1}-\widetilde{v}^{[k]}(t_{k+1})\|_{\mathbb{H}}\big]$.
In view of \eqref{ddd}, it is clear that
\begin{align*}
	\mathbb{E}\big[\left\|u_{k+1}-S_{\Delta t}u_k\right\|^2_{\mathbb{H}}\big]\leq C\Delta t
\end{align*}
and
\begin{align*}
	v_{k+1}
	&=v_k+\Delta tS^{N-k}_{\Delta t}F(S_{\Delta t}u_k)+J_k\Delta W_k\\
	&\quad+\Delta t^2S^{N-k}_{\Delta t}DF(S_{\Delta t}u_k)S_{\Delta t}F(u_{k+1})\\
	&\quad+\Delta tS^{N-k}_{\Delta t}\int_{0}^{1}(1-\theta)D^2F(u_k^{\theta})(u_{k+1}-S_{\Delta t}u_k,u_{k+1}-S_{\Delta t}u_k)\dd \theta,
\end{align*}
which leads to
 \begin{equation}\label{syyu}
 	\begin{split}
 	\mathbb{E}\big[\|v_{k+1}-\widetilde{v}^{[k]}(t_{k+1})\|_{\mathbb{H}}\big]&\leq\Delta t^2\left\|S_{\Delta t}\right\|^{N-k+1}_{\mathcal{L}(\mathbb{H};\mathbb{H})}\|DF\|_{\mathbf{C}_b}\mathbb{E}\big[\left\|F(u_{k+1})\right\|_{\mathbb{H}}\big]\\
 		&\quad+\Delta t\left\|S_{\Delta t}\right\|^{N-k}_{\mathcal{L}(\mathbb{H};\mathbb{H})}\|D^2F\|_{\mathbf{C}_b}\mathbb{E}\big[\left\|u_{k+1}-S_{\Delta t}u_k\right\|_{\mathbb{H}}^2\big]\\
 		&\leq C\Delta t^2\exp\Big(-\frac{g\sigma_0}{2}(N-k)\Delta t\Big),
 	\end{split}
 \end{equation}
 where $u_k^{\theta}:=S_{\Delta t}u_k+\theta(u_{k+1}-S_{\Delta t}u_k)$ .
 
 \textbf{(1b). Estimate of }$\mathbb{E}\big[\|\widetilde{v}^{[k]}(t)-v_k\|_{\mathbb{H}}^{2i}\big],i=1,2.$ 
 In view of the Burkholder--Davis--Gundy-type inequality, we have
	\begin{equation}\label{wwes}
		\begin{split}
			\mathbb{E}\big[\|\widetilde{v}^{[k]}(t)-v_k\|_{\mathbb{H}}^{2i}\big]
			&\leq C\mathbb{E}\Big[\Big(\int_{t_k}^{t}\left\|S_{\Delta t}^{N-k}F(S_{\Delta t}u_k)\right\|_{\mathbb{H}}\dd s\Big)^{2i}\Big]\\
			&\quad+C\mathbb{E}\Big[\Big(\int_{t_k}^{t}\left\|J_k\right\|_{HS(U_0;\mathbb{H})}^{2}\dd s\Big)^i\Big]\\
			&\leq C\Delta t^i\exp\big(-ig\sigma_0(N-k)\Delta t\big)\qquad\forall\, t\in[t_k,t_{k+1}].
		\end{split}
	\end{equation}
	
	 \textbf{(1c). Estimate of }$\mathbb{E}\big[\|\widetilde{v}^{[k]}(t)\|_{\mathbb{H}}^2\big]$. It follows from the It\^o isometry that
	 \begin{align*}
	 	\mathbb{E}\big[\|\widetilde{v}^{[k]}(t)\|_{\mathbb{H}}^2\big]\leq C\exp\big(-g\sigma_0(N-k)\Delta t\big)\qquad\forall\,t\in[t_k,t_{k+1}].
	 \end{align*}

\textbf{Step 2.}
Since $V(0,v_N)=\varphi(u_N)$ and
		$V(t_N,v_0)=U(t_N,S(-t_N)S^N_{\Delta t}u_0)$,
	we split the weak error into the following terms
	\begin{equation}\label{weak_error}
		\begin{split}
		\mathbb{E}\big[\varphi(u_N)\big]-\mathbb{E}\big[\varphi(u(t_N))\big]
		&	=\underbrace{\mathbb{E}\big[V(t_N,v_0)\big]-\mathbb{E}\big[U(t_N,u_0)\big]}_{=:\mathcal{H}}\\	&\quad+\sum_{k=0}^{N-1}\Big\{\underbrace{\mathbb{E}\big[V(t_N-t_{k+1},v_{k+1})\big]-\mathbb{E}\big[V(t_N-t_{k+1},\widetilde{v}^{[k]}(t_{k+1}))\big]}_{=:\mathcal{J}^{(k)}}\Big\}\\
		&\quad+\sum_{k=0}^{N-1}\Big\{\underbrace{\mathbb{E}\big[V(t_N-t_{k+1},\widetilde{v}^{[k]}(t_{k+1}))\big]-\mathbb{E}\big[V(t_N-t_k,v_k)\big]}_{=:\mathcal{A}^{(k)}}\Big\}.
			\end{split}
	\end{equation}
	
	For the term $\mathcal{H}$, we have
	\begin{equation}\label{4}
		\begin{split}
		\big|\mathcal{H}\big|
		&=\Big|\mathbb{E}\Big[\int_{0}^{1}DU\big(t_N,u_0+\theta \big(S(-t_N)S_{\Delta t}^Nu_0-u_0\big)\big)S(-t_N)\big(S_{\Delta t}^N-S(t_N)\big)u_0\dd \theta\Big]\Big|\\
		&\leq C\exp\big(K_1^{(1)}N\Delta t\big)\|S_{\Delta t}^N-S(t_N)\|_{\mathcal{L}(\mathcal{D}(M^2);\mathbb{H})}\mathbb{E}\big[\|u_0\|_{\mathcal{D}(M^2)}\big]\\
		&\leq C\Delta t \exp\Big(-\frac{\widetilde{g}\sigma_0-2K_1^{(1)}}{2}N\Delta t\Big).
	\end{split}
	\end{equation}
	
Notice that
\begin{align*}
	\big|\mathcal{J}^{(k)}\big|
	&\leq C\exp\big(K_1^{(1)}(N-k-1)\Delta t\big)\mathbb{E}\big[\|v_{k+1}-\widetilde{v}^{[k]}(t_{k+1})\|_{\mathbb{H}}\big]\\
	&\leq C\Delta t^2\exp\Big(-\frac{g\sigma_0-2K_1^{(1)}}{2}(N-k-1)\Delta t\Big),
\end{align*}
which implies that
\begin{align}\label{3}
	\Big| \sum_{k=0}^{N-1}\mathcal{J}^{(k)}\Big|\leq C\Delta t^2 \sum_{k=0}^{N-1}\exp\Big(-\frac{g\sigma_0-2K_1^{(1)}}{2}(N-k-1)\Delta t\Big)\leq C\Delta t.
\end{align}

	 In view of \eqref{spp1}, we now use the It\^o formula on $V(t_N-t,\widetilde{v}^{[k]}(t))$ to obtain
  that for $t\in[t_k,t_{k+1}]$,
	\begin{align*}
		\dd V(t_N-t,\widetilde{v}^{[k]}(t))=\mathcal{A}^{(k)}_1(t)\dd t+\mathcal{A}^{(k)}_2(t)\dd t
		+\langle DV(t_N-t,\widetilde{v}^{[k]}(t)),J_k\dd W(t)\rangle,
	\end{align*}
where
	\begin{align*}
		\mathcal{A}^{(k)}_1(t)&:=-\frac{1}{2}\text{Tr}\Big[D^2V(t_N-t,\widetilde{v}^{[k]}(t))\Lambda(t_N-t,\widetilde{v}^{[k]}(t))\Lambda(t_N-t,\widetilde{v}^{[k]}(t))^*\Big]\\
		&\quad+\frac{1}{2}\text{Tr}\Big[D^2V(t_N-t,\widetilde{v}^{[k]}(t))\big(J_kQ^{\frac{1}{2}}\big)\big(J_kQ^{\frac{1}{2}}\big)^*\Big],\\
		\mathcal{A}^{(k)}_2(t)&:= DV(t_N-t,\widetilde{v}^{[k]}(t))\Big(S_{\Delta t}^{N-k}F(S_{\Delta t}u_k)-\Upsilon(t_N-t,\widetilde{v}^{[k]}(t))\Big).
	\end{align*}
	
	Let $\widetilde{J}_k:=\Delta tS_{\Delta t}^{N-k}DF(S_{\Delta t}u_k))S_{\Delta t}B(u_k)$. Further, we decompose $\mathcal{A}^{(k)}_1(t)$ and $\mathcal{A}^{(k)}_2(t)$, respectively, as
	\begin{align*}
	\mathcal{A}^{(k)}_1(t)
		&=\mathcal{A}^{(k)}_{1,1}(t)+\mathcal{A}^{(k)}_{1,2}(t)+\mathcal{A}^{(k)}_{1,3}(t)+\mathcal{A}_{1,4}^{(k)}(t)\quad\text{and}\quad
	\mathcal{A}^{(k)}_2(t)=	\mathcal{A}^{(k)}_{2,1}(t)+	\mathcal{A}_{2,2}^{(k)}(t)+\mathcal{A}_{2,3}^{(k)}(t),
\end{align*}
	where 
	\begin{align*}
			\mathcal{A}^{(k)}_{1,1}(t)&:=-\frac{1}{2}\text{Tr}\Big[D^2V(t_N-t,v_k)\Lambda(t_N-t,\widetilde{v}^{[k]}(t))\Lambda(t_N-t,\widetilde{v}^{[k]}(t))^*\Big]\\
			&\quad+\frac{1}{2}\text{Tr}\Big[D^2V(t_N-t,v_k)\Lambda(t_N-t,v_k)\Lambda(t_N-t,v_k)^*\Big],\\
			\mathcal{A}^{(k)}_{1,2}(t)&:=\frac{1}{2}\text{Tr}\Big[D^2V(t_N-t,v_k)\Lambda(t_N-t,\widetilde{v}^{[k]}(t))\Lambda(t_N-t,\widetilde{v}^{[k]}(t))^*\Big]\\
			&\quad-\frac{1}{2}\text{Tr}\Big[D^2V(t_N-t,\widetilde{v}^{[k]}(t))\Lambda(t_N-t,\widetilde{v}^{[k]}(t))\Lambda(t_N-t,\widetilde{v}^{[k]}(t))^*\Big]\\
			&\quad+\frac{1}{2}\text{Tr}\Big[D^2V(t_N-t,\widetilde{v}^{[k]}(t))\Lambda(t_N-t,v_k)\Lambda(t_N-t,v_k)^*\Big]\\
			&\quad-\frac{1}{2}\text{Tr}\Big[D^2V(t_N-t,v_k)\Lambda(t_N-t,v_k)\Lambda(t_N-t,v_k)^*\Big],\\		
	\mathcal{A}^{(k)}_{1,3}(t)&:=-\frac{1}{2}\text{Tr}\Big[D^2V(t_N-t,\widetilde{v}^{[k]}(t))\Lambda(t_N-t,v_k)\Lambda(t_N-t,v_k)^*\Big]\\
		&\quad+\frac{1}{2}\text{Tr}\Big[D^2V(t_N-t,\widetilde{v}^{[k]}(t))\big(S_{\Delta t}^{N-k}B(u_k)Q^{\frac{1}{2}}\big)\big(S_{\Delta t}^{N-k}B(u_k)Q^{\frac{1}{2}}\big)^*\Big],\\
		\mathcal{A}^{(k)}_{1,4}(t)&:=\frac{1}{2}\text{Tr}\Big[D^2V(t_N-t,\widetilde{v}^{[k]}(t))\big(\widetilde{J}_kQ^{\frac{1}{2}}\big)\big(J_kQ^{\frac{1}{2}}\big)^*\Big]\\
		&\quad+\frac{1}{2}\text{Tr}\Big[D^2V(t_N-t,\widetilde{v}^{[k]}(t))\big(S_{\Delta t}^{N-k}B(u_k)Q^{\frac{1}{2}}\big)\big(\widetilde{J}_kQ^{\frac{1}{2}}\big)^*\Big],\\
		\mathcal{A}^{(k)}_{2,1}(t)&:= DV(t_N-t,\widetilde{v}^{[k]}(t))\Big(S_{\Delta t}^{N-k}F(S_{\Delta t}u_k)-\Upsilon(t_N-t,v_k)\Big),\\
		\mathcal{A}_{2,2}^{(k)}(t)&:= DV(t_N-t,v_k)\Big(\Upsilon(t_N-t,v_k)-\Upsilon(t_N-t,\widetilde{v}^{[k]}(t))\Big),\\
		\mathcal{A}_{2,3}^{(k)}(t)&:=\Big(DV(t_N-t,\widetilde{v}^{[k]}(t))-DV(t_N-t,v_k)\Big)\Big(\Upsilon(t_N-t,v_k)-\Upsilon(t_N-t,\widetilde{v}^{[k]}(t))\Big).
	\end{align*}

\textbf{Step 3. Estimate of } $\mathcal{A}^{(k)}_1(t)$. 

	\textbf{(3a).} 
	For fixed $r\in[t_k,t_{k+1}]$, let $G^{r,i}(\widetilde{v}^{[k]}(t)):=D^2V(t_N-r,v_k)\left(\Lambda(t_N-r,\widetilde{v}^{[k]}(t))e_i,\Lambda(t_N-r,\widetilde{v}^{[k]}(t))e_i\right).$
Notice that for every $h\in\mathbb{H}$,
\begin{align*}
	DG^{r,i}(\widetilde{v}^{[k]}(t))h&=2D^2V(t_N-r,v_k)\Big(\Lambda^{(1)}(t_N-r,\widetilde{v}^{[k]}(t),h)e_i,\Lambda(t_N-r,\widetilde{v}^{[k]}(t))e_i\Big),
\end{align*}
where $\Lambda^{(1)}(t_N-r,\widetilde{v}^{[k]}(t),h):=S(t_N-r)DB(S(-t_N+r)\widetilde{v}^{[k]}(t))S(-t_N+r)hQ^{\frac{1}{2}}.$
Besides, for every $h_1,h_2\in\mathbb{H}$, we have
\begin{align*}
	D^2G^{r,i}(\widetilde{v}^{[k]}(t))(h_1,h_2)&=2D^2V(t_N-r,v_k)\Big(\Lambda^{(1)}(t_N-r,\widetilde{v}^{[k]}(t),h)e_i,\Lambda^{(1)}(t_N-r,\widetilde{v}^{[k]}(t),h_2)e_i\Big)\\
&\quad+2D^2V(t_N-r,v_k)\Big(\Lambda^{(2)}(t_N-r,\widetilde{v}^{[k]}(t),h_1,h_2)e_i,\Lambda(t_N-r,\widetilde{v}^{[k]}(t))e_i\Big),
\end{align*}
where $\Lambda^{(2)}(t_N-r,\widetilde{v}^{[k]}(t),h_1,h_2):=S(t_N-r)D^2B(S(-t_N+r)\widetilde{v}^{[k]}(t))\big(S(-t_N+r)h_1,S(-t_N+r)h_2\big)Q^{\frac{1}{2}}.$

	We use the It\^o formula on $G^{r,i}(\widetilde{v}^{[k]}(t))$  and take the expectation to obtain that
	\begin{align*}
	\mathbb{E}\big[\mathcal{A}_{1,1}^{(k)}(t)\big]
		&=-\frac{1}{2}\sum_{i=1}^{\infty}\Big\{
		\int_{t_k}^{t}\mathbb{E}\big[DG^{t,i}(\widetilde{v}^{[k]}(s))S_{\Delta t}^{N-k}F(S_{\Delta t}u_k)\big]\dd s\\
		&\quad+\frac{1}{2}\int_{t_k}^{t}\sum_{j=1}^{\infty}\mathbb{E}\big[D^2G^{t,i}(\widetilde{v}^{[k]}(s))\big(J_kQ^{\frac{1}{2}}e_j,J_kQ^{\frac{1}{2}}e_j\big)\big]\dd s\Big\},
	\end{align*} 
	which leads to
	\begin{equation}\label{A11}
		\begin{split}
	\big|\mathbb{E}\big[\mathcal{A}_{1,1}^{(k)}(t)\big]\big|&\leq C\Delta t\exp\Big(-\big((g-1)\sigma_0-K^{(2)}_1\big)(N-k)\Delta t\Big)\qquad\forall \,t\in[t_k,t_{k+1}].
\end{split}
	\end{equation}

\textbf{(3b).}	For the term $\mathcal{A}^{(k)}_{1,2}(t)$, notice that
	\begin{align*}
		 \mathcal{A}_{1,2}^{(k)}(t)
		 &=-\frac{1}{2}\sum_{i=1}^{\infty}\bigg[\int_{0}^{1}D^3V\big(t_N-t,v_k+\theta(\widetilde{v}^{[k]}(t)-v_k)\big)\\
		 & \qquad\qquad\qquad\Big(\Lambda(t_N-t,\widetilde{v}^{[k]}(t))e_i,\Lambda(t_N-t,\widetilde{v}^{[k]}(t))e_i,\widetilde{v}^{[k]}(t)-v_k\Big)\dd\theta\bigg]\\
		 &\quad+\frac{1}{2}\sum_{i=1}^{\infty}\bigg[\int_{0}^{1}D^3V\big(t_N-t,v_k+\theta(\widetilde{v}^{[k]}(t)-v_k)\big)\\
		 &\qquad\qquad\qquad\Big(\Lambda(t_N-t,v_k)e_i,\Lambda(t_N-t,v_k)e_i,\widetilde{v}^{[k]}(t)-v_k\Big)\dd\theta\bigg]\\
		 	&=\frac{1}{2}\sum_{i=1}^{\infty}\bigg[\int_{0}^{1}D^3V\big(t_N-t,v_k+\theta(\widetilde{v}^{[k]}(t)-v_k)\big)\\
		 &\qquad\qquad\qquad\Big(\Lambda(t_N-t,v_k)e_i-\Lambda(t_N-t,\widetilde{v}^{[k]}(t))e_i,\Lambda(t_N-t,\widetilde{v}^{[k]}(t))e_i,\widetilde{v}^{[k]}(t)-v_k\Big)\dd\theta\bigg]\\
		&\quad+\frac{1}{2}\sum_{i=1}^{\infty}\bigg[\int_{0}^{1}D^3V\big(t_N-t,v_k+\theta(\widetilde{v}^{[k]}(t)-v_k)\big)\\
		& \qquad\qquad\qquad\Big(\Lambda(t_N-t,v_k)e_i, \Lambda(t_N-t,v_k)e_i-\Lambda(t_N-t,\widetilde{v}^{[k]}(t))e_i,\widetilde{v}^{[k]}(t)-v_k\Big)\dd\theta\bigg].
	\end{align*}
	Hence, it holds that
	\begin{align}\label{A12}
		\mathbb{E}\big[\big|\mathcal{A}_{1,2}^{(k)}(t)\big|\big]\leq  C\Delta t\exp\Big(-\frac{(3g-4)\sigma_0-2K^{(4)}}{2}(N-k)\Delta t\Big)\qquad\forall\,t\in[t_k,t_{k+1}].
	\end{align}

%In view of \eqref{111} and \eqref{112}, we have
%\begin{align}\label{11}
%		\Big|\sum_{k=0}^{N-1}\mathbb{E}\Big[\int_{t_k}^{t_{k+1}}\mathcal{A}_{1,1}^{(k)}(t)\dd t\Big]\Big|\leq  C\Delta t.
%\end{align}

	\textbf{(3c).} Let $z_k(t):=S_{\Delta t}^{N-k}B(u_k)-S(t_N-t)B(S(-t_N+t)v_k)$ where $t\in[t_k,t_{k+1}].$
	Notice that
	\begin{align*}
		z_k(t)
		&=S_{\Delta t}^{N-k}\big(B(u_k)-B(S(-t_N+t)v_k)\big)+\big(S_{\Delta t}^{N-k}- S(t_N-t)\big)B(S(-t_N+t)v_k),
	\end{align*}
	which implies that
	\begin{align*}
		&\quad\mathbb{E}\big[\|z_k(t)\|^2_{HS(U_0;\mathbb{H})}\big]\\
		&\leq 2\|S_{\Delta t}\|_{\mathcal{L}(\mathbb{H};\mathbb{H})}^{2N-2k}\mathbb{E}\big[\|B(u_k)-B(S(-t_N+t)S_{\Delta t}^{N-k}u_k)\|^2_{HS(U_0;\mathbb{H})}\big]\\
		&\quad+2\|S_{\Delta t}^{N-k}- S(t_N-t)\|^2 _{\mathcal{L}(\mathcal{D}(M^2);\mathbb{H})}\mathbb{E}\big[\|B(S(-t_N+t)v_k)\|^2_{HS(U_0;\mathcal{D}(M^2))}\big]\\
		&\leq C\|S_{\Delta t}\|_{\mathcal{L}(\mathbb{H};\mathbb{H})}^{2N-2k}\|S(-t_N+t)\|^2_{\mathcal{L}(\mathbb{H};\mathbb{H})}\|S(t_N-t)-S_{\Delta t}^{N-k}\|^2_{\mathcal{L}(\mathcal{D}(M^2);\mathbb{H})}\mathbb{E}\big[\|u_k\|^2_{\mathcal{D}(M^2)}\big]\\
		&\quad+2\|S_{\Delta t}^{N-k}- S(t_N-t)\|^2_{\mathcal{L}(\mathcal{D}(M^2);\mathbb{H})}\mathbb{E}\big[\|B(S(-t_N+t)v_k)\|^2_{HS(U_0;\mathcal{D}(M^2))}\big]\\
		&\leq C\Delta t^2\exp\Big(-(\widetilde{g}+g-2)\sigma_0(N-k)\Delta t\Big).
	\end{align*}
	
	Hence, we have
	\begin{align*}
		\mathcal{A}^{(k)}_{1,3}(t)
		&=\frac{1}{2}\text{Tr}\Big[D^2V(t_N-t,\widetilde{v}^{[k]}(t))\big(z_k(t)Q^{\frac{1}{2}}\big)\Lambda(t_N-t,v_k)^*\Big]\\
		&\quad+\frac{1}{2}\text{Tr}\Big[D^2V(t_N-t,\widetilde{v}^{[k]}(t))\big(S_{\Delta t}^{N-k}B(u_k)Q^{\frac{1}{2}}\big)\big(z_k(t)Q^{\frac{1}{2}}\big)^*\Big],
	\end{align*}
	which leads to 
	\begin{align}\label{A13}
		\mathbb{E}\big[\big|\mathcal{A}_{1,3}^{(k)}(t)\big|\big]
		&\leq C\Delta t \exp\Big(-\frac{(\widetilde{g}+2g-4)\sigma_0-2K^{(2)}_1}{2}(N-k)\Delta t\Big).
	\end{align}

	\textbf{(3d).}	For the term $\mathcal{A}^{(k)}_{1,4}(t)$, we have
	\begin{align}\label{A14}
		\mathbb{E}\big[\big|\mathcal{A}_{1,4}^{(k)}(t)\big|\big]&\leq C\exp\Big(\big(\sigma_0+K_1^{(2)}\big)(N-k)\Delta t\Big)\Big(\mathbb{E}\big[\|\widetilde{J}_k\|_{HS(U_0;\mathbb{H})}^2\big]\Big)^{\frac{1}{2}}\notag\\
		&\qquad\cdot\Big(\mathbb{E}\big[\|S_{\Delta t}^{N-k}B(u_k)+J_k\|_{HS(U_0;\mathbb{H})}^2\big]\Big)^{\frac{1}{2}}\\
		&\leq C\Delta t\exp\Big(-\big((g-1)\sigma_0-K^{(2)}_1\big)(N-k)\Delta t\Big)\qquad\forall \,t\in[t_k,t_{k+1}].\notag
	\end{align}
	
	In view of \eqref{A11}-\eqref{A14}, we arrive at
		\begin{align}\label{1}
		\left|\sum_{k=0}^{N-1}\mathbb{E}\left[\int_{t_k}^{t_{k+1}}\mathcal{A}_{1}^{(k)}(t)\dd t\right]\right|\leq C\Delta t.
	\end{align}

	\textbf{Step 4. Estimate of }$\mathcal{A}^{(k)}_2(t)$.
	
\textbf{(4a).}	Denote $\widetilde{z}_k(t):=S_{\Delta t}^{N-k}F(S_{\Delta t}u_k)-\Upsilon(t_N-t,v_k)$ for $t\in[t_k,t_{k+1}].$
	Notice that
	\begin{align*}
		\widetilde{z}_k(t)
		&=S_{\Delta t}^{N-k}\Big(F(S_{\Delta t}u_k)-F(S(-t_N+t)v_k)\Big)+\Big(S_{\Delta t}^{N-k}-S(t_N-t)\Big)F(S(-t_N+t)v_k),
	\end{align*}
	which implies that
	\begin{align*}
		\mathbb{E}\left[\|\widetilde{z}_k(t)\|^2_{\mathbb{H}}\right]
		&\leq C\|S_{\Delta t}\|_{\mathcal{L}(\mathbb{H};\mathbb{H})}^{2N-2k}\mathbb{E}\left[\|S_{\Delta t}u_k-S(-t_N+t)S_{\Delta t}^{N-k}u_k\|^2_{\mathbb{H}}\right]\\
		&\quad+C\|S_{\Delta t}^{N-k}- S(t_N-t)\|^2_{\mathcal{L}(\mathcal{D}(M^2);\mathbb{H})} \mathbb{E}\left[\|F(S(-t_N+t)v_k)\|^2_{\mathcal{D}(M^2)}\right]\\
		&\leq C\|S_{\Delta t}\|_{\mathcal{L}(\mathbb{H};\mathbb{H})}^{2N-2k}\|S(-t_N+t)\|^2_{\mathcal{L}(\mathbb{H};\mathbb{H})}\|S(t_N-t)-S_{\Delta t}^{N-k-1}\|^2_{\mathcal{L}(\mathcal{D}(M^2);\mathbb{H})}\\
		&\qquad\qquad\qquad\qquad\qquad\cdot\|S_{\Delta t}\|^2_{\mathcal{L}(\mathcal{D}(M^2);\mathcal{D}(M^2))}\mathbb{E}\left[\|u_k\|^2_{\mathcal{D}(M^2)}\right]\\
		&\quad+C\|S_{\Delta t}^{N-k}- S(t_N-t)\|^2_{\mathcal{L}(\mathcal{D}(M^2);\mathbb{H})}\mathbb{E}\left[\|F(S(-t_N+t)v_k)\|^2_{\mathcal{D}(M^2)}\right]\\
		&\leq C\Delta t^2\exp\Big(-(\widetilde{g}+g-2)(N-k)\sigma_0\Delta t\Big).
	\end{align*}
	Therefore, it holds that
	\begin{equation}\label{A21}
		\begin{split}
		\mathbb{E}\big[\big|\mathcal{A}^{(k)}_{2,1}(t)\big|\big]&\leq C\exp\Big(K_1^{(1)}(N-k)\Delta t\Big)\Big(\mathbb{E}\big[\|\widetilde{z}_k(t)\|_{\mathbb{H}}^2\big]\Big)^{\frac{1}{2}}\\
		&\leq C\Delta t\exp\Big(-\frac{(\widetilde{g}+g-2)\sigma_0-2K_1^{(1)}}{2}(N-k)\Delta t\Big)\qquad\forall\, t\in[t_k,t_{k+1}].
			\end{split}
	\end{equation}

\textbf{(4b).}	
For fixed $r\in[t_k,t_{k+1}]$, let $\widetilde{G}^r(\widetilde{v}^{[k]}(t)):=DV(t_N-r,v_k)\Upsilon(t_N-r,\widetilde{v}^{[k]}(t)),t\in[t_k,t_{k+1}].$
Notice that
$$
	D\widetilde{G}^r(\widetilde{v}^{[k]}(t))h=DV(t_N-r,v_k)D\Upsilon(t_N-r,\widetilde{v}^{[k]}(t))h\qquad\forall \,h\in\mathbb{H},
$$
$$
	D^2\widetilde{G}^r(\widetilde{v}^{[k]}(t))(h_1,h_2)=DV(t_N-r,v_k)D^2\Upsilon(t_N-r,\widetilde{v}^{[k]}(t))(h_1,h_2)\qquad\forall \,h_1,h_2\in\mathbb{H},
$$
where 
$$
	D\Upsilon(t_N-r,\widetilde{v}^{[k]}(t))h=S(t_N-r)DF(S(-t_N+r)\widetilde{v}^{[k]}(t))S(-t_N+r)h,
$$
$$
	D^2\Upsilon(t_N-r,\widetilde{v}^{[k]}(t))(h_1,h_2)=S(t_N-r)D^2F(S(-t_N+r)\widetilde{v}^{[k]}(t))\big(S(-t_N+r)h_1,S(-t_N+r)h_2\big).
$$

		By the It\^o formula on $\widetilde{G}^{r}(\widetilde{v}^{[k]}(t))$  and taking the expectation, it yields that
\begin{align*}
	\mathbb{E}\big[\mathcal{A}_{2,2}^{(k)}(t)\big]
	&=-\int_{t_k}^{t}\mathbb{E}\big[D\widetilde{G}^t(\widetilde{v}^{[k]}(s))S_{\Delta t}^{N-k}F(S_{\Delta t}u_k)\big]\dd s\\
	&\quad+\frac{1}{2}\int_{t_k}^{t}\sum_{i=1}^{\infty}\mathbb{E}\big[D^2\widetilde{G}^t(\widetilde{v}^{[k]}(s))\big(J_kQ^{\frac{1}{2}}e_i,J_kQ^{\frac{1}{2}}e_i\big)\big]\dd s,
\end{align*}
by which we obtain
\begin{align}\label{A22}
	\big|\mathbb{E}\big[\mathcal{A}_{2,2}^{(k)}(t)\big]\big|&\leq C\Delta t\exp\big(K_1^{(1)}(N-k)\Delta t\big)\|S_{\Delta t}\|^{N-k}_{\mathcal{L}(\mathbb{H};\mathbb{H})}\Big(\mathbb{E}\big[\|F(S_{\Delta t}u_k)\|_{\mathbb{H}}\big]\notag\\
	&\qquad+\|S(-t_N+t)\|_{\mathcal{L}(\mathbb{H};\mathbb{H})}\|S_{\Delta t}\|^{N-k}_{\mathcal{L}(\mathbb{H};\mathbb{H})}\mathbb{E}\big[\|B(u_k)\|^2_{HS(U_0;\mathbb{H})}\big]\Big)\\
	&\leq C\Delta t\exp\Big(-\big((g-1)\sigma_0-K_1^{(1)}\big)(N-k)\Delta t\Big)\qquad \forall\,t\in[t_k,t_{k+1}].\notag
\end{align}

\textbf{(4c).}	For the term $\mathcal{A}_{2,3}^{(k)}(t)$, we have
\begin{align*}
	\mathcal{A}_{2,3}^{(k)}(t)&=\int_{0}^{1}(1-\theta)D^2V(t_N-t,v_k+\theta(\widetilde{v}^{[k]}(t)-v_k))\\
	&\qquad\quad\big(\Upsilon(t_N-t,v_k)-\Upsilon(t_N-t,\widetilde{v}^{[k]}(t)),\widetilde{v}^{[k]}(t)-v_k\big)\dd\theta,
\end{align*}
which leads to 
\begin{equation}\label{A23}
	\begin{split}
	\mathbb{E}\big[\big|\mathcal{A}_{2,3}^{(k)}(t)\big|\big]&\leq C\exp\Big((\sigma_0+K_1^{(2)})(N-k)\Delta t)\Big)\mathbb{E}\big[\|\widetilde{v}^{[k]}(t)-v_k\|^2_{\mathbb{H}}\big]\\
	&\leq C\Delta t\exp\Big(-\big((g-1)\sigma_0-K^{(2)}_1\big)(N-k)\Delta t\Big)\qquad \forall\,t\in[t_k,t_{k+1}].
\end{split}
\end{equation}
In view of \eqref{A21}-\eqref{A23}, we conclude that
\begin{align}\label{2}	\left|\sum_{k=0}^{N-1}\mathbb{E}\left[\int_{t_k}^{t_{k+1}}\mathcal{A}_{2}^{(k)}(t)\dd t\right]\right|\leq C\Delta t.
\end{align}

	Combining \textbf{Step 2}, \textbf{Step 3}, \textbf{Step 4}, we finish the proof.
\end{proof}

\begin{rema}\label{djk}
	(i) For the case of additive noise,  it is clear that $\mathcal{A}_{1,1}^{(k)}(\cdot)\equiv0$ and $\mathcal{A}_{1,2}^{(k)}(\cdot)\equiv0$. Hence, we can obtain the weak convergence order one for every $\varphi\in\mathbf{C}_b^2(\mathbb{H};\mathbb{R})$.\\
	(ii) Under the condition of \cref{zx}, we can also prove the strong convergence order for \eqref{semi-numerical}, i.e., $\sup_{n\geq 0}\big(\mathbb{E}\big[\|u(t_n)-u_n\|_{\mathbb{H}}^2\big]\big)^{\frac{1}{2}}\leq C\Delta t^{\frac{1}{2}}$, which is half of the weak convergence order.
\end{rema}

\subsection{Applications}
In this part, we give some applications of the weak convergence result, including the error estimate between the invariant measures $\pi$ and $\pi^{\Delta t}$, the strong law of large numbers and central limit theorem related to time average of the numerical solution, and the convergence order of multi-level Monte Carlo estimator.
\subsubsection{Convergence order between invariant measures $\pi$ and $\pi^{\Delta t}$}
Since it is a nontrivial task to find the exact invariant measure $\pi$, an alternative way is to approximate it using the numerical invariant measure $\pi^{\Delta t}$. In the following, we investigate the convergence order of the numerical invariant measure to the exact one.

\begin{prop}
	 Under the conditions of \cref{zx}, we have 
	\begin{align*}
		\Big|\int_{\mathbb{H}}\varphi(y)\dd\pi(y)-\int_{\mathbb{H}}\varphi(y)\dd\pi^{\Delta t}\Big|\leq C\Delta t\qquad\forall\,\varphi\in\mathbf{C}^3_b(\mathbb{H};\mathbb{R}).
	\end{align*}
\end{prop}
\begin{proof}
	Since $\{u(t)\}_{t\geq 0}$ and $\{u_k\}_{k\in\mathbb{N}}$ are ergodic, we have
	\begin{align*}
		\lim_{T\rightarrow\infty}\frac{1}{T}\int_{0}^{T}\mathbb{E}\big[\varphi(u(t))\big]\dd t&=\int_{\mathbb{H}}\varphi(y)\dd\pi(y),\\
		\lim_{M\rightarrow\infty}\frac{1}{M}\sum_{k=0}^{M-1}\mathbb{E}\big[\varphi(u_k)\big]&=\int_{\mathbb{H}}\varphi(y)\dd\pi^{\Delta t}(y)
	\end{align*}
for any $\varphi\in\mathbf{C}^3_b(\mathbb{H};\mathbb{R})$. Hence, it yields that for sufficiently small $\Delta t>0$,
\begin{align*}
	\Big|\int_{\mathbb{H}}\varphi(y)\dd\pi(y)-\int_{\mathbb{H}}\varphi(y)\dd\pi^{\Delta t}(y)\Big|\leq\lim_{\substack{M\rightarrow\infty,\\T=M\Delta t\rightarrow\infty}}\frac{1}{T}\sum_{k=0}^{M-1}\int_{t_k}^{t_{k+1}}\big|\mathbb{E}\big[\varphi(u(t))\big]-\mathbb{E}\big[\varphi(u_k)\big]\big|\dd t\leq C\Delta t.
\end{align*}
\end{proof}

\subsubsection{Limit theorems related to the time average of the numerical solution}
Based on the convergence order of the semi-implicit Euler scheme, we study the limit theorems of the time average $\frac{1}{k}\sum_{i=0}^{k-1}f(u_k)$ of the numerical solution, including the strong law of large numbers  and the central limit theorem.

In view of \cref{smzpp}, \cref{dfz21}, Remark \ref{djk}(ii), and \cite[Theorem 3.1]{chen2023probabilistic}, we first present the strong law of large numbers for the time average in the following proposition.
\begin{prop}
	Let $\widetilde{q}\in\mathbb{N}_+,p\geq 1$, and $\gamma\in(0,1]$ with $p+\gamma\leq \widetilde{q}$.
	Under the conditions of  \cref{dfz21} with $q=\widetilde{q}$, and  \cref{zx}, we have
	\begin{align*}
		\lim_{\Delta t\rightarrow0}\lim_{k\rightarrow\infty}\frac{1}{k}\sum_{i=0}^{k-1}f(u_i)\stackrel{a.s.}{=}\pi(f)\qquad\forall\, f\in \mathcal{C}_{p,\gamma}(\mathbb{H};\mathbb{R}).
	\end{align*}
\end{prop}

Moreover, based on  \cref{smzpp}, \cref{dfz21}, \cref{zx}, Remark \ref{djk}(ii), and \cite[Theorem 3.2 and Remark 2]{chen2023probabilistic}, we also obtain the following central limit theorem. Recall that $v^2=2\pi((f-\pi(f))\int_{0}^{\infty}(P_tf-\pi(f))\dd t)$ is given in \cref{drr}.
\begin{prop}
	(i) Let $q_1\in\mathbb{N}_+,r_1\geq 2$,$\gamma_1\in(0,1]$, $k_1=\lceil\Delta t^{-1-2\lambda_1}\rceil$ with $\lambda_1\in(0,\frac{1}{2}\gamma_1)$, and $p_1\geq 1$ with
	$
	\big(3\vee\frac{1+\lambda_1}{\lambda_1}\big)(p_1+3\gamma_1)\leq (2q_1)\wedge r_1.$
	Under the conditions of  \cref{smzpp} with $r=r_1$,  \cref{dfz21} with $q=q_1$, and \cref{zx}, we have that for any
	$f\in \mathcal{C}_{p_1,\gamma_1}(\mathbb{H};\mathbb{R})$, 
	\begin{align*}
		\frac{1}{\sqrt{k_1\Delta t}}\sum_{i=0}^{k_1-1}\big(f(u_i)-\pi(f)\big)\Delta t\stackrel{d}{\rightarrow}\mathcal{N}(0,v^2)\quad\text{as }\Delta t\rightarrow0.
	\end{align*}\\
	(ii) Let $q_2\in\mathbb{N}_+,r_2\geq 2$,$\gamma_2\in(0,1]$,  $k_2=\lceil\Delta t^{-1-2\lambda_2}\rceil$ with $\lambda_2\in(0,\gamma_2)$, and $p_2\geq 1$ with
$
		\big(3\vee\frac{1+\lambda_2}{\lambda_2}\big)(p_2+3\gamma_2)\leq (2q_2)\wedge r_2.$
Under the conditions of  \cref{smzpp} with $r=r_2$,  \cref{dfz21} with $q=q_2$, and \cref{zx}, we have that for any
 $f\in \mathbf{C}_b^3(\mathbb{H};\mathbb{R})$, 
\begin{align*}
	\frac{1}{\sqrt{k_2\Delta t}}\sum_{i=0}^{k_2-1}\big(f(u_i)-\pi(f)\big)\Delta t\stackrel{d}{\rightarrow}\mathcal{N}(0,v^2)\quad\text{as }\Delta t\rightarrow0.
\end{align*}
\end{prop}

\begin{rema}
	Notice that for $f\in \mathbf{C}_b^3(\mathbb{H};\mathbb{R})$, the solution of the Poisson equation $\widetilde{\mathcal{L}}\varphi=f-\pi(f)$ is $\varphi=-\int_{0}^{\infty}\big(P_tf-\pi(f)\big)\dd t$, where  the generator $\widetilde{\mathcal{L}}$ of \eqref{www} is defined as
	\begin{align*}
		\widetilde{\mathcal{L}}f(v)=\langle Df(v),Mv-\sigma_0v+F(v)\rangle+\frac{1}{2}\text{Tr}\big[D^2f(v)\big(B(v)Q^{\frac{1}{2}}\big)\big(B(v)Q^{\frac{1}{2}}\big)^{*}\big]\qquad\forall\, v\in\mathcal{D}(M).
	\end{align*}
	 Since $\varphi\widetilde{\mathcal{L}}\varphi=\frac{1}{2}\widetilde{\mathcal{L}}\varphi^2-\frac{1}{2}\|D\varphi B\|^2_{HS(U_0;\mathbb{R})}$ and $\pi(\widetilde{\mathcal{L}}\varphi^2)=0$, we have $v^2=-2\pi(\varphi\widetilde{\mathcal{L}}\varphi)=\pi(\|D\varphi B\|^2_{HS(U_0;\mathbb{R})})$.
\end{rema}

\subsubsection{Convergence order of multi-level Monte Carlo estimator}
For fixed $T>0$, it is difficult to obtain the exact value of $\mathbb{E}\big[\varphi(u(T))\big]$, which is usually numerically approximated by the multi-level Monte Carlo method. Below, we give the convergence order of the multi-level Monte Carlo estimator with the help of the convergence order of \eqref{semi-numerical}.

% To this end, we give the following theorem concerning the strong convergence order of \eqref{semi-numerical}.
 
% \begin{assumption}\label{assu}
% 	Assume that 
% 	\begin{align*}
% 		\|F(h_1)-F(h_2)\|_{\mathbb{H}}\leq C_F\|h_1-h_2\|_{\mathbb{H}}\qquad&\forall\, h_1,h_2\in\mathbb{H},\\
% 		\|F(h)\|_{\mathcal{D}(M)}\leq \widetilde{C}_F(1+\|h\|_{\mathcal{D}(M)})\qquad&\forall \,h\in\mathcal{D}(M).
% 	\end{align*}
% \end{assumption}
 
% \begin{theo}\label{strong order}
% 	Let Assumption \ref{ass} with $0\leq \rho\leq 2$ and Assumption \ref{assu} hold. Let  $u_0\in L^2(\Omega;\mathcal{D}(M^2))$ and $\sigma_0>\max\{\widetilde{\alpha}_{F,2}+\widetilde{\alpha}^2_{B,2},\alpha_{F,0}+2\alpha^2_{B,0}+2\}$. For sufficiently small $\Delta t>0$, there exists a positive constant $C$ independent of $\Delta t$ such that
% 	\begin{align*}
% 		\sup_{n\geq 0}\mathbb{E}\big[\|u(t_n)-u_n\|_{\mathbb{H}}^2\big]\leq C\Delta t.
% 	\end{align*}
% \end{theo}
% \begin{proof}
% 	The proof is similar to that of \cite[Theorem 4.3]{CHJ2019}, thus we omit it here.
% \end{proof}

 Consider  a sequence of equidistant nested time partitions $\Delta =\{\Delta ^l\}_{l\in\mathbb{N}}$ where $\Delta ^l=\{t_n^l=nT\Delta t_l,n=0,\cdots,N_l\}$, $N_l=2^l$, and $\Delta t_l=2^{-l}$. Then we get a sequence of approximations $\{u_T^l:=u_{N_l}\}_{l\in\mathbb{N}}$ for $u(T)$.
In view of \cref{zx} and Remark \ref{djk}(ii), there exists a $L_0\in\mathbb{N}$ such that for $\varphi\in\mathbf{C}_b^3(\mathbb{H};\mathbb{R})$,
\begin{align}\label{sdt1}
	\big|\mathbb{E}\big[\varphi(u(T))\big]-\mathbb{E}\big[\varphi(u_T^{L_0+l})\big]\big|\leq C\Delta t_{L_0+l}\qquad\forall\,l\geq 0,
\end{align} 
and
\begin{equation}\label{sdt2}
	\begin{split}
	\text{Var}\big[\varphi(u_T^{L_0+l})-\varphi(u_{T}^{L_0+l-1})\big]&\leq \big\|\varphi(u_T^{L_0+l})-\varphi(u_T^{L_0+l-1})\big\|^2_{L^2(\Omega;\mathbb{R})}\\
	&\leq 2\big\|\varphi(u_T^{L_0+l})-\varphi(u(T))\big\|^2_{L^2(\Omega;\mathbb{R})}\\
	&\quad+2\big\|\varphi(u_T^{L_0+l-1})-\varphi(u(T))\big\|^2_{L^2(\Omega;\mathbb{R})}\\
	&\leq C\Delta t_{L_0+l}+C\Delta t_{L_0+l-1}\leq C\Delta t_{L_0+l}\qquad\forall\,l\geq 1.
		\end{split}
\end{equation}
% We observe that  \cref{zx} and Remark \ref{djk}(ii)  imply that in the setting of \cite[Theorem 1]{Lang} that $Y=\varphi(u(T))$ and there exists a $L_0\in\mathbb{N}$ such that $Y_l=\varphi(u_T^{L_0+l}),a_l=2^{-L_0-l}$ and $\eta=\frac{1}{2}$.

 For a given $L\in\mathbb{N}_+$, we define the multi-level Monte Carlo estimator of $\mathbb{E}\big[\varphi(u_T^{L_0+l})\big]$ by
 \begin{align*}
 	E^L\big[\varphi(u_T^{L_0+l})\big]:=E_{M_0}\big[\varphi(u_T^{L_0})\big]+\sum_{l=1}^{L}E_{M_l}\big[\varphi(u_T^{L_0+l})-\varphi(u_T^{L_0+l-1})\big],
 \end{align*}
where $E_{M}[X]$ is the Monte Carlo estimator of $X$.

At this stage, we can obtain the convergence order of multi-level Monte Carlo estimator in the following corollary.

\begin{coro}
Under the condition of \cref{zx}, the multilevel Monte Carlo estimator on level $L>0$ satisfies for any $\epsilon>0$ at the final time $T$,
	\begin{align}\label{dwwq}
		\big\|\mathbb{E}\big[\varphi(u(T))\big]-E^L\big[\varphi(u_T^{L_0+L})\big]\big\|_{L^2(\Omega;\mathbb{R})}\leq (C+C(1+\zeta(1+\epsilon))^{\frac{1}{2}})\Delta t_{L_0+L}\qquad\forall\,\varphi\in\mathbf{C}_b^3(\mathbb{H};\mathbb{R}),
	\end{align}
with sample sizes given by $M_l=\lceil 2^{L_0+2L-l}l^{1+\epsilon}\rceil$, $l=1,\cdots,L$, and $M_0= 2^{2L_0+2L}$,where $\lceil\cdot\rceil$ denotes the rounding to the next large integer and $\zeta$ the Riemann zeta function.
\end{coro}
\begin{proof}
In view of \cite[Lemma 2]{Lang}, we use \eqref{sdt1} and \eqref{sdt2} to obtain that
\begin{align*}	&\quad\big\|\mathbb{E}\big[\varphi(u(T))\big]-E^L\big[\varphi(u_T^{L_0+L})\big]\big\|_{L^2(\Omega;\mathbb{R})}\\
	&\leq \big|\mathbb{E}\big[\varphi(u(T))\big]-\mathbb{E}\big[\varphi(u_T^{L_0+L})\big]\big|+\Big(M_0^{-1}\text{Var}[u_T^{L_0}]+\sum_{l=1}^{L}M_l^{-1}\text{Var}\big[\varphi(u_T^{L_0+l})-\varphi(u_{T}^{L_0+l-1})\big]\Big)^{\frac{1}{2}}\\
	&\leq C\Delta t_{L_0+L}+\Big(C\Delta t^2_{L_0+L}+\sum_{l=1}^{L}l^{-(1+\epsilon)}\Delta t_{L_0+2L-l}\Delta t_{L_0+l}\Big)^{\frac{1}{2}}\\
	&=C\Delta t_{L_0+L}+C\Big(1+\sum_{l=1}^{L}l^{-(1+\epsilon)}\Big)^{\frac{1}{2}}\Delta t_{L_0+L}.
\end{align*}
 Thus, we finish the proof.
\end{proof}

\appendix
\section{The Proof of Lemma \ref{sfu4}}
Let $v_k:=S_{\Delta t}^kv_0$ and $v(t):=S(t)v_0$. Then
\begin{align}\label{sw4}
	v_k=v_{k-1}-\Delta t\sigma_0v_k+\Delta tMv_k
\end{align}
and 
\begin{align}
	v(t_k)=v(t_{k-1})+\int_{t_{k-1}}^{t_k}(M-\sigma_0)v(s)\dd s.	
\end{align}

(1)	
\textbf{Case 1.} $i=0$:  By applying $\langle\cdot,v_k\rangle_{\mathbb{H}}$ on both sides of \eqref{sw4}, we have
\begin{align*}
	\frac{1}{2}\Big(\|v_k\|^2_{\mathbb{H}}-\|v_{k-1}\|^2_{\mathbb{H}}\Big)+\frac{1}{2}\|v_k-v_{k-1}\|^2_{\mathbb{H}}=-\Delta t\sigma_0\|v_k\|_{\mathbb{H}}^2,
\end{align*}
which implies that 
\begin{align}\label{z141}
	\|v_k\|^2_{\mathbb{H}}\leq\frac{1}{1+2\Delta t\sigma_0}\|v_{k-1}\|^2_{\mathbb{H}}\leq\cdots\leq\frac{1}{\left(1+2\Delta t\sigma_0\right)^k}\|v_0\|_{\mathbb{H}}^2.
\end{align}
There exists a constant $\Delta t_g>0$ such that
\begin{align*}
	\frac{1}{1+2\sigma_0\Delta t}\leq e^{-g\sigma_0\Delta t}\qquad\forall\,\Delta t\in(0,\Delta t_g).
\end{align*} 
Hence, we have 
\begin{align}\label{z14}
	\|v_k\|^2_{\mathbb{H}}\leq e^{-gk\sigma_0\Delta t}\|v_0\|_{\mathbb{H}}^2,
\end{align}
i.e., 
$$\|S^k_{\Delta t}\|_{\mathcal{L}(\mathbb{H};\mathbb{H})}\leq e^{-\frac{g}{2}k\sigma_0\Delta t}.$$

\textbf{Case 2.} $i\in\{1,2\}$: 	
In view of \eqref{sw4}, one gets
\begin{align}\label{xc04}
	M^iv_k=M^iv_{k-1}-\Delta t\sigma_0M^iv_k+\Delta tM^{i+1}v_k.
\end{align}
By applying $\langle\cdot,M^iv_k\rangle_{\mathbb{H}}$ on both sides of \eqref{xc04}, it yields that
\begin{align*}
	\frac{1}{2}\Big(\|M^iv_k\|^2_{\mathbb{H}}-\|M^iv_{k-1}\|^2_{\mathbb{H}}\Big)+\frac{1}{2}\|M^iv_k-M^iv_{k-1}\|^2_{\mathbb{H}}=-\Delta t\sigma_0\|M^iv_k\|^2_{\mathbb{H}},
\end{align*}
which leads to
\begin{align}\label{z214}
	\|M^iv_k\|^2_{\mathbb{H}}\leq\frac{1}{(1+2\sigma_0\Delta t)^k}\|M^iv_0\|^2_{\mathbb{H}}\leq e^{-gk\sigma_0\Delta t}\|M^iv_0\|^2_{\mathbb{H}}.
\end{align}

Combining \eqref{z14} and \eqref{z214}, it holds that $\|v_k\|_{\mathcal{D}(M^i)}\leq e^{-\frac{g}{2}k\sigma_0\Delta t} \|v_0\|_{\mathcal{D}(M^i)},$
i.e., $\|S^k_{\Delta t}\|_{\mathcal{L}(\mathcal{D}(M^i);\mathcal{D}(M^i))}\leq e^{-\frac{g}{2}k\sigma_0\Delta t}$.

%\textbf{Case 3.} $i=2$:
%In view of \eqref{sw4}, we have
%\begin{align}\label{xc4}
%	M^2v_k=M^2v_{k-1}-\sigma_0\Delta tM^2v_k+\Delta tM^3v_k.
%\end{align}
%By applying $\langle\cdot,M^2v_k\rangle_{\mathbb{H}}$ on both sides of \eqref{xc4}, we have
%\begin{align*}
%	\frac{1}{2}\Big(\|M^2v_k\|^2_{\mathbb{H}}-\|M^2v_{k-1}\|^2_{\mathbb{H}}\Big)+\frac{1}{2}\|M^2v_k-M^2v_{k-1}\|^2_{\mathbb{H}}=-\sigma_0\Delta t\|M^2v_k\|^2_{\mathbb{H}},
%\end{align*}
%which implies that
%\begin{align}\label{z24}
%	\|M^2v_k\|^2_{\mathbb{H}}\leq\frac{1}{1+2\sigma_0\Delta t}\|M^2v_{k-1}\|^2_{\mathbb{H}}\leq\cdots\leq\frac{1}{(1+2\sigma_0\Delta t)^k}\|M^2v_0\|^2_{\mathbb{H}}\leq e^{-gk\sigma_0\Delta t}\|M^2v_0\|^2_{\mathbb{H}}.
%\end{align}
%
%Combining \eqref{z14} and \eqref{z24}, we have
%\begin{align*}
%\|v_k\|_{\mathcal{D}(M^2)}\leq e^{-\frac{g}{2}k\sigma_0\Delta t} \|v_0\|_{\mathcal{D}(M^2)},
%\end{align*}
%which implies that 
%\begin{align*}
%	\|S_{\Delta t}^k\|_{\mathcal{L}(\mathcal{D}(M^2);\mathcal{D}(M^2))}\leq e^{-\frac{g}{2}k\sigma_0\Delta t}.
%\end{align*}

(2) Notice that
\begin{align*}
	\|(e^{-t\sigma_0I}-I)v\|_{\mathbb{H}}= (1-e^{-t\sigma_0})\|v\|_{\mathbb{H}}\leq \sigma_0 t\|v\|_{\mathbb{H}}\leq \sigma_0t\|v\|_{\mathcal{D}(M)}\qquad\forall\,v\in\mathcal{D}(M),
\end{align*}
i.e., $\|e^{-t\sigma_0I}-I\|_{\mathcal{L}(\mathcal{D}(M);\mathbb{H})}\leq \sigma_0t$.
 
Therefore, 
 \begin{align*}
 	\|S(t)-e^{tM}\|_{\mathcal{L}(\mathcal{D}(M);\mathbb{H})}
 	\leq \|e^{tM}\|_{\mathcal{L}(\mathbb{H};\mathbb{H})}\|e^{-t\sigma_0I}-I\|_{\mathcal{L}(\mathcal{D}(M);\mathbb{H})}\leq \sigma_0t,
 \end{align*}
by which it holds that
\begin{align*}
	\|S(t)-I\|_{\mathcal{L}(\mathcal{D}(M);\mathbb{H})}\leq \|S(t)-e^{tM}\|_{\mathcal{L}(\mathcal{D}(M);\mathbb{H})}+\|e^{tM}-I\|_{\mathcal{L}(\mathcal{D}(M);\mathbb{H})}\leq Ct\qquad\forall \,t\in[0,\infty).
\end{align*}

(3) \textbf{Step 1.} For any $v_0\in\mathcal{D}(M^2)$, denote
$e_k:=v(t_k)-v_k=(S(t_k)-S_{\Delta t}^k)v_0, k\in\mathbb{N}.$
It is clear that $\|e_0\|_{\mathbb{H}}=0$, i.e., 
\begin{align}\label{A9}
	\|S_{\Delta t}^0-S(t_0)\|_{\mathcal{L}(\mathcal{D}(M^2);\mathbb{H})}=0.
\end{align}

For $k\in\mathbb{N}_+$, we have
\begin{equation}\label{mmmm4}
	\begin{split}
		e_k&=v(t_{k-1})-v_{k-1}+\int_{t_{k-1}}^{t_k}\Big(M-\sigma_0I\Big)v(s)\dd s+\Delta t\sigma_0v_k-\Delta tMv_k\\
		&=e_{k-1}+\int_{t_{k-1}}^{t_k}\Big(M-\sigma_0I\Big)v(s)\dd s-\Delta t \Big(M-\sigma_0I\Big)v_k\\
		&=e_{k-1}+\Delta t\Big(M-\sigma_0I\Big)e_k+\int_{t_{k-1}}^{t_k}\Big(M-\sigma_0I\Big)\Big(v(s)-v(t_k)\Big)\dd s.
	\end{split}
\end{equation}

By applying $\langle \cdot,e_k\rangle_{\mathbb{H}}$ on both sides of \eqref{mmmm4}, one gets
\begin{align*}
	&\quad\frac{1}{2}\Big(\|e_k\|^2_{\mathbb{H}}-\|e_{k-1}\|^2_{\mathbb{H}}\Big)+\frac{1}{2}\|e_k-e_{k-1}\|^2_{\mathbb{H}}\\
	&=-\sigma_0\Delta t\|e_k\|^2_{\mathbb{H}}+\int_{t_{k-1}}^{t_k}\left\langle \big(M-\sigma_0I\big)\big(v(s)-v(t_k)\big),e_k\right\rangle_{\mathbb{H}}\dd s\\
	&=-\sigma_0\Delta t\|e_k\|^2_{\mathbb{H}}-\int_{t_{k-1}}^{t_k}\langle v(s)-v(t_k),(M+\sigma_0I)e_k\rangle_{\mathbb{H}}\dd s\\
	&=-\sigma_0\Delta t\|e_k\|^2_{\mathbb{H}}+\int_{t_{k-1}}^{t_k}\int_{s}^{t_k}\left\langle (M-\sigma_0I)v(r),(M+\sigma_0I)e_k\right\rangle_{\mathbb{H}}\dd r\dd s\\
	&\leq -\sigma_0\Delta t\|e_k\|^2_{\mathbb{H}}+\int_{t_{k-1}}^{t_k}\int_{s}^{t_k}\Big( 2\sigma_0\|Mv(r)\|_{\mathbb{H}}+\|M^2v(r)\|_{\mathbb{H}}+\sigma_0^2\|v(r)\|_{\mathbb{H}}\Big)\|e_k\|_{\mathbb{H}}\dd r\dd s\\
	&\leq-\sigma_0\Delta t\|e_k\|^2_{\mathbb{H}}+ C(\sigma_0+1)^2\Delta t^2\|e_k\|_{\mathbb{H}}\sup_{t_{k-1}\leq t<t_{k}}\|v(t)\|_{\mathcal{D}(M^2)}\\
	&\leq -\sigma_0\Delta t\|e_k\|^2_{\mathbb{H}}+C(\sigma_0+1)^2\Delta t^2\|e_k\|_{\mathbb{H}}\|v_0\|_{\mathcal{D}(M^2)}e^{-\sigma_0(k-1)\Delta t}\\
	&\leq -\frac{1}{2}\sigma_0\Delta t\|e_k\|^2_{\mathbb{H}}+Ce^{-2\sigma_0(k-1)\Delta t}\Delta t^3\|v_0\|^2_{\mathcal{D}(M^2)}.
\end{align*}

Notice that there exists a positive constant $\Delta t_{\widetilde{g}}$ such that
\begin{align*}
	\frac{1}{1+\sigma_0\Delta t}\leq e^{-\widetilde{g}\sigma_0\Delta t}\qquad\forall\,\Delta t\in(0,\Delta t_{\widetilde{g}}),
\end{align*}
which leads to
\begin{align*}
	\|e_k\|^2_{\mathbb{H}}&\leq \frac{1}{1+\sigma_0\Delta t}\|e_{k-1}\|^2_{\mathbb{H}}+\frac{C\Delta t^3}{1+\sigma_0\Delta t}e^{-2\sigma_0(k-1)\Delta t}\|v_0\|^2_{\mathcal{D}(M^2)}\\
	&\leq e^{-\widetilde{g}\sigma_0\Delta t}\|e_{k-1}\|^2_{\mathbb{H}}+C\Delta t^3e^{-2\sigma_0\Delta t(k-1)-\widetilde{g}\sigma_0\Delta t}\|v_0\|^2_{\mathcal{D}(M^2)}\\
	&\leq\cdots
	\leq C\Delta t^3\|v_0\|^2_{\mathcal{D}(M^2)}e^{-2k\sigma_0\Delta t}\sum_{i=1}^{k}e^{(2-\widetilde{g})i\sigma_0\Delta t}\\
	&=C\Delta t^3\|v_0\|^2_{\mathcal{D}(M^2)}e^{-\widetilde{g}k\sigma_0\Delta t}\sum_{i=1}^{k}e^{-(2-\widetilde{g})(k-i)\sigma_0\Delta t}\\
	&\leq C\Delta t^2\|v_0\|^2_{\mathcal{D}(M^2)}e^{-\widetilde{g}k\sigma_0\Delta t}\qquad\forall~\Delta t\in (0,\Delta t_{\widetilde{g}}).
\end{align*}
Combining with \eqref{A9}, it yields that for $k\in\mathbb{N}$,
\begin{align*}
	\|S_{\Delta t}^{k}-S(t_k)\|_{\mathcal{L}(D(M^2;\mathbb{H}))}\leq Ce^{-\frac{\widetilde{g}}{2}k\sigma_0\Delta t}\Delta t\qquad\forall~\Delta t\in (0,\Delta t_{\widetilde{g}}).
\end{align*}

\textbf{Step 2}.
Notice that
\begin{align*}
	\|(S_{\Delta t}^0-S(t))v_0\|_{\mathbb{H}}
	&\leq\|S(t)-I\|_{\mathcal{L}(\mathcal{D}(M);\mathbb{H})}\|v_0\|_{\mathcal{D}(M)}\leq C\Delta t\|v_0\|_{\mathcal{D}(M^2)}\qquad \forall\,t\in[t_0,t_1],
\end{align*}
and
\begin{align*}
	\|(S_{\Delta t}^k-S(t))v_0\|_{\mathbb{H}}&\leq\|(S_{\Delta t}^k-S(t_k))v_0\|_{\mathbb{H}}+\|(S(t_k)-S(t))v_0\|_{\mathbb{H}}\\
%	&=\|(S_{\Delta t}^k-S(t_k))v_0\|_{\mathbb{H}}+\|S(t_k\wedge t)(I-S(|t-t_k|))v_0\|_{\mathbb{H}}\\
	&\leq \|S_{\Delta t}^k-S(t_k)\|_{\mathcal{L}(\mathcal{D}(M^2);\mathbb{H})}\|v_0\|_{\mathcal{D}(M^2)}\\
	&\quad+\|S(t_k\wedge  t)\|_{\mathcal{L}(\mathbb{H};\mathbb{H})}\|I-S(|t-t_k|)\|_{\mathcal{L}(\mathcal{D}(M);\mathbb{H})}\|v_0\|_{\mathcal{D}(M)}\\
	&\leq Ce^{-\frac{\widetilde{g}}{2}k\sigma_0\Delta t}  \Delta t\|v_0\|_{\mathcal{D}(M^2)}\qquad\forall \,t\in[t_{k-1},t_{k+1}],\,k\in\mathbb{N}_+.
\end{align*}

\bibliographystyle{IEEEtran}	
\bibliography{references}

\end{document}